\documentclass[12pt,reqno, tbtags]{amsart}
\usepackage{mathrsfs}
\usepackage{amsmath}
\usepackage{epsfig}
\usepackage{amsfonts}
\usepackage{amssymb}
\usepackage{amsthm}
\usepackage{amscd}
\usepackage[all]{xy}
\usepackage[T1]{fontenc}
\usepackage{rotating}
\usepackage{lscape}
\usepackage{fancyhdr}
\usepackage{amsbsy}
\usepackage{verbatim}
\usepackage{moreverb}
\usepackage{rotating}

\newcommand{\ra}{\rightarrow}

\newcommand{\xra}{\xrightarrow}

\renewcommand{\S}{{\rm S}}

\newcommand{\R}{{\mathbb R}}
\newcommand{\sm}{\wedge}

\newcommand{\id}{{\rm id}}

\newcommand{\Z}{{\mathbb Z}}

\newcommand{\C}{{\mathbb{C}}}

\newcommand{\OO}{{\mathcal{O}}}

\newcommand{\HH}{{\mathcal{H}}}
\newcommand{\GL}{GL}

\newcommand{\E}{{\rm E}}
\newcommand{\LL}{{\mathcal{L}}}
\newcommand{\PPic}{{\mathcal{P}\emph{ic}}}
\newcommand{\Sp}{{\mathit{Sp}}}
\newcommand{\CC}{{\mathcal{C}}}
\newcommand{\RR}{{\mathcal{R}}}
\newcommand{\MM}{{\mathcal{M}}}

\newcommand{\pt}{\ast}
\newcommand{\A}{{\textbf A}}
\newcommand{\tw}{{\textrm{tw}}}
\newcommand{\Ho}{{\textrm{Ho}}}
\newcommand{\sCat}{{\mathit{sCat}}}
\newcommand{\infCat}{\infty\textrm{-\textit{Cat}}}
\newcommand{\sSet}{{\mathit{sSet}}}
\newcommand{\qCat}{{\mathit{qCat}}}
\newcommand{\hGpd}{{\mathit{hGpd}}}

\newcommand{\TT}{{\mathbb{T}}}
\newcommand{\Th}{\textrm{Th}}
\newcommand{\Rmod}{\textrm{R-mod}}
\newcommand{\Smod}{\mathrm{S^0}\textrm{-mod}}
\newcommand{\res}{\textrm{res}}

\renewcommand{\SS}{S^0}
\renewcommand{\AA}{{\textrm{A}}}

\newcommand{\HZ}{\rm{H}\mathbb{Z}}

\newcommand{\nn}{\nonumber}
\newcommand{\nid}{\noindent}

\newcommand{\twiddles}{\begin{center} $\sim\sim\sim\sim\sim$ \end{center}}

\renewcommand{\:}{\! :}

\newtheorem{theorem}{Theorem}[section]
\newtheorem{prop}[theorem]{Proposition}
\newtheorem{lemma}[theorem]{Lemma}
\newtheorem{cor}[theorem]{Corollary}

\theoremstyle{definition}
\newtheorem{defn}[theorem]{Definition}
\newtheorem{rechar}[theorem]{Recharacterization}
\newtheorem{summary}[theorem]{Summary}
\newtheorem{note}[theorem]{Note}

\theoremstyle{remark}
\newtheorem{remark}[theorem]{Remark}
\newtheorem{example}[theorem]{Example}
\newtheorem{conj}[theorem]{Conjecture}

\DeclareMathOperator*{\colim}{colim}
\DeclareMathOperator*{\holim}{holim}
\DeclareMathOperator*{\hocolim}{hocolim}
\DeclareMathOperator{\Aut}{Aut}

\DeclareMathOperator{\Pic}{Pic}
\DeclareMathOperator{\Haunt}{Haunt}
\DeclareMathOperator{\Specter}{Specter}
\DeclareMathOperator{\Hom}{Hom}
\DeclareMathOperator{\haut}{haut}

\CompileMatrices

\begin{document}

\title{Twisted Parametrized Stable Homotopy Theory}
\author{Christopher L. Douglas}
\date{12 June 2005}
\address{Department of Mathematics, M. I. T., Cambridge, MA
02139}
\email{cdouglas@math.mit.edu}
\thanks{The author was partially supported by a Clay Mathematics Institute Liftoff Fellowship.}

\vspace*{-12pt}

\begin{abstract}
We introduce a framework, twisted parametrized stable homotopy theory, for
describing semi-infinite homotopy types.  A twisted parametrized
spectrum is a section of a bundle whose fibre is the category of spectra.
We define these bundles in terms of modules over a stack of parametrized
spectra and in terms of diagrams of simplicial categories.  We present
a classification of bundles of categories of spectra and of the associated
twisted parametrized spectra.  Though twisted parametrized spectra do
not have global homotopy types and therefore do not have generalized homology invariants
in the usual sense, they do admit generalized-homology-type invariants for certain
commutative ring spectra.  We describe this invariant theory and in particular
note that under mild hypotheses, a twisted parametrized spectrum will have
cyclically graded homology invariants and will also have K-theory and complex
bordism invariants, as expected for a Floer or semi-infinite homotopy type.
We discuss the association of a twist of parametrized homotopy theory
to a polarized infinite-dimensional manifold and present a conjectural,
explicit realization of this twist in terms of parametrized semi-infinitely
indexed spectra.
\end{abstract}

\vspace*{-12pt}

\maketitle

\vspace*{-12pt}

\setcounter{tocdepth}{3}
\tableofcontents

\vspace*{-18pt}

\section{Introduction}  

\subsection{Background and Motivation}

Despite its widespread use and compelling application to problems in symplectic topology and gauge theory, Floer homology remains rather a mystery.  The very existence of the Floer homology of an infinite-dimensional manifold depends on delicate and haphazard properties of the manifold and of the flow associated to a Floer function; moreover, it is completely unknown how the Floer homology depends on the choice of Floer function.  Confronted with this situation, Cohen, Jones, and Segal~\cite{cjs} asked whether it would be possible to build a Floer homotopy type encoding the relevant data from the manifold and the function in such a way that a homology functor would recover Floer homology.  Besides elucidating the structure of Floer theory and clarifying its dependence on the Floer function, such a Floer homotopy type would immediately provide other invariants such as Floer $K$-theory and Floer bordism.  Cohen, Jones, and Segal suggested that prospectra might encode some of the Floer data; though this thought proved useful, in retrospect it is clear that prospectra can only account for the Floer homotopy types of trivially polarized manifolds---this restriction on the polarization partially accounts for the difficulty Cohen, Jones, and Segal had finding examples of Floer prospectra.

The purpose of this paper is not to answer the Floer homotopy question, but to introduce a framework, namely twisted parametrized stable homotopy theory, that is a necessarily component of any description of Floer or semi-infinite homotopy.  That some twisted form of homotopy theory was needed to account for nontrivial polarizations was first realized by Furuta~\cite{furuta}, and it will turn out that the twisted
space he wrote down (as a conjectural model for the Seiberg-Witten-Floer homotopy type of $T^3$) is a very specialized example of our twisted parametrized spectra.  A twisted parametrized spectrum is a section of a bundle whose fibre is the category of spectra, and as such it has the same relationship to an ordinary parametrized spectrum as a section of a line bundle has to a function.  A polarized infinite-dimensional manifold has a naturally associated bundle of categories of spectra, and the fundamental ansatz is that geometric information about such a manifold (its semi-infinite homotopy type, for example) involves this bundle and its sections.  In this paper, we present the theory of twisted parametrized spectra, including various definitions, characterizations, and classifications of them, a thorough description of their homotopy-theoretic invariants, and an overview of their relationship to infinite-dimensional polarized manifolds.  The specific association of a twisted parametrized homotopy type to a manifold with Floer function is the subject of ongoing work with Mike Hopkins and Ciprian Manolescu~\cite{douglashopkins, douglasmanolescu}.

Applications to semi-infinite homotopy theory aside, twisted parametrized
spectra provide a natural and informative generalization of ordinary parametrized spectra.  As the study of modular functions begs for a theory of modular forms, so the study of parametrized spectra is intrinsically accompanied by a theory of twisted parametrized spectra.  Besides expanding
the range of homotopy types that can be described in parametrized topology,
the twisted theory provides new and refined perspectives on such classical topics as Thom spectra on loop spaces and parametrized homotopy and cohomotopy invariants.  Indeed, many of the natural algebraic invariants arising in twisted parametrized homotopy theory are novel even when applied to ordinary parametrized spectra.

I would like to thank especially Mike Hopkins for insightful and inspiring questions and indispensable pointers, Bill Dwyer for fruitful suggestions and encouraging words, and Jacob Lurie for technical help and much headache-saving advice.

\subsection{Overview}

A twisted parametrized spectrum is a section of a bundle whose fibre is the category of spectra.  There are two natural (and in the end equivalent)
definitions of such bundles; the first definition is based on the notion of invertible sheaves of categories of spectra and the second definition formalizes the concept of local systems of categories of spectra.
These two approaches to twisted parametrized homotopy theory are developed, respectively, in sections~\ref{sectionsheaves} and~\ref{sectionlocsys}.

An invertible sheaf on an algebraic variety $X$ over a ring $R$ is a rank-one locally free module over the structure sheaf of $R$-valued functions on $X$.  The basic trope of this paper is to replace the ring $R$ by the category $\Sp$ of spectra.  A $\Sp$-valued function is naturally interpreted as a parametrized spectrum, and an invertible sheaf of categories of spectra is therefore a rank-one locally free module over the structure stack of parametrized spectra on $X$.  This structure stack of parametrized spectra is described in section~\ref{ssps} and modules over the structure stack are discussed in section~\ref{modstack}.  Invertible sheaves of categories of spectra are referred to as \emph{haunts}.  Section~\ref{haunts} presents various examples of haunts and motivates the fundamental classification result: haunts on a space $X$ are classified by homotopy classes of maps from $X$ to a deloop $B\Pic(\SS)$ of the realization of the Picard category of the sphere spectrum.

The global sections of a haunt are the \emph{twisted parametrized spectra} or \emph{specters} for short.  Section~\ref{specters} presents numerous examples of specters, including a thorough homotopy-theoretic description of the specters on low-dimensional spheres, and outlines the general classification: there is a bundle \mbox{$E\Pic(\SS) \times_{\Pic(\SS)} \Sp_w$} over $B\Pic(\SS)$ whose fibre is the realization $\Sp_w$ of the weak equivalence subcategory of the category of spectra, and specters for a fixed haunt are classified by homotopy classes of lifts of the classifying map $X \ra B\Pic(\SS)$ of the haunt to the total space of the bundle \mbox{$E\Pic(\SS) \times_{\Pic(\SS)} \Sp_w$}.  The invertible sheaves approach to twisted parametrized stable homotopy theory is summarized in table~\ref{tableframework}.

A local system is a sheaf $L$ with the property that the ring of sections $L(U)$ over any sufficiently small open set $U$ is isomorphic to a fixed ring $R$.  By analogy, a local system of categories of spectra is a stack $M$ whose category of sections $M(U)$ over any small open set $U$ is appropriately equivalent to the category $\Sp$ of spectra.  In order to describe the relevant notion of equivalence, we need a category of categories of spectra---such a category is typically referred to as a `homotopy theory of homotopy theories'. The category of simplicial categories provides a convenient framework and is the subject of section~\ref{htpyhtpy}.  In section~\ref{diagsimp} we formally define local systems of categories of spectra in terms of diagrams of simplicial categories.  It turns out that a local system of categories of spectra is equivalent to an invertible sheaf of categories of spectra; a haunt may therefore be freely interpreted as either a local system or an invertible sheaf.  The bulk of section~\ref{haunts2} is devoted to establishing, in the context of diagrams of simplicial categories, the classification of haunts announced in section~\ref{haunts}.  Correspondingly, section~\ref{specters2} proves a categorical classification result for specters which has as a corollary the classification discussed in section~\ref{specters}.  The local systems approach to twisted parametrized stable homotopy theory is summarized in table~\ref{secondtableframework}.

Building on the classification results of section~\ref{diagsimp}, section~\ref{sectionthom} relates haunts and specters to more traditional homotopy-theoretic concepts.  We describe in section~\ref{automspec} the homotopy type of the Picard category of the sphere spectrum and thereby identify the classifying space for haunts as a deloop of the classifying space for stable spherical fibrations.  In section~\ref{specmod} we use this identification to recharacterize haunts in terms of $\AA_{\infty}$ ring spectra arising as Thom spectra on loop spaces; we describe the associated categories of specters as module categories for the ring Thom spectra.

With the basic theory in place, we devote section~\ref{sectioninvariants} to a description of the invariants of specters.  Unlike ordinary parametrized spectra, specters do not have global homotopy or cohomotopy types.
Nevertheless, given a commutative ring spectrum $R$ there is a parametrized generalized homology functor that associates to a specter a corresponding "$R$-specter"---this association is described in section~\ref{Rhaunt}.  For appropriate choices of the ring spectrum $R$, the $R$-specter \emph{will} have a global homotopy type and we can therefore define the $R$-homology groups of the original specter.  Section~\ref{Rhomology} presents various examples of
these specter invariants and describes a spectral sequence for their computation.

The final section of the paper connects twisted parametrized homotopy theory to the geometry of polarized Hilbert manifolds.  In section~\ref{polhilb} we delineate the relevant notions of polarization, namely real and complex polarizations of real and complex Hilbert spaces, and we describe the homotopy types of the corresponding structure groups and classifying spaces.  The classifying space for real polarizations of a real Hilbert space maps to the classifying space for categories of spectra.  A real polarized Hilbert bundle therefore has an associated haunt; section~\ref{symppol} describes this association and presents examples of naturally occurring polarized manifolds and their haunts.
When a real polarization of a real Hilbert bundle lifts to a real polarization of a complex Hilbert bundle, the corresponding haunt and its specters have a greatly simplified invariant theory.  For example, in section~\ref{unitpol} we note that the homology of a specter for such a polarization is necessarily graded by a finite cyclic group---this explains the idea that semi-infinite and Floer homology theories are cyclically, not integrally, graded.  Moreover we prove that under mild hypotheses, specters for such a polarization (and therefore a large class of `semi-infinite homotopy types') have global $\HZ$-, $K$-, and $MU$-homology invariants.  This explains and substantially generalizes the Cohen-Jones-Segal~\cite{cjs} observation that trivially polarized manifolds can admit associated Floer homology, $K$-theory, and complex bordism invariants.

In the concluding section~\ref{siis} we define the notion of a spectrum indexed not on the finite-dimensional subspaces of an infinite-dimensional vector space but on the infinite-dimensional subspaces of a Hilbert space that are compatible with a fixed polarization.  Though they are conjecturally equivalent to ordinary spectra, the resulting semi-infinitely indexed spectra are explicitly and naturally linked to the geometry of infinite-dimensional manifolds.  In particular, we expect that the category of parametrized semi-infinitely indexed spectra on a polarized Hilbert bundle is a model for the category of specters for the haunt associated to the polarized bundle.

\section{Invertible Sheaves of Categories of Spectra}  \label{sectionsheaves} 

Let $X$ be an algebraic variety over a ring $R$.  The structure sheaf $\OO_X$ can be described as the sheaf of $R$-valued functions on $X$.  The most fundamental $\OO_X$-modules are the invertible or locally free rank-one modules.  These modules, which we will often think of as line bundles, are classified by the first cohomology group of $X$ with coefficients in the sheaf $\OO_X^{\times}$ of invertible functions.  Denote by $\LL(c)$ the line bundle associated to the cohomology class $c \in H^1(X;\OO_X^{\times})$.  A global section of $\LL(c)$ is determined by a 0-cochain on $X$ whose coboundary is $c$, that is by an element $f \in C^0(X;\OO_X)$ such that $\delta f = c$.

Twisted parametrized stable homotopy theory is a precise analog of these algebraic concepts: the ring $R$ is replaced by the category $\Sp$ of spectra, a "categorical semi-ring" under the wedge and smash products.  The set of $\Sp$-valued functions on a space $X$ is naturally interpreted as the category of parametrized spectra on $X$, and the structure "sheaf" is therefore the structure stack $\OO_X$ of parametrized spectra.  There is a notion of locally free rank-one module over the structure stack of parametrized spectra and we refer to such modules, briefly, as \emph{haunts}.  Haunts are classified by the first cohomology group of $X$ with coefficients in the so-called Picard stack $\PPic(\SS)$ of invertible parametrized spectra.  The fundamental objects of twisted parametrized stable homotopy theory are the global sections of a haunt; these global sections are the twisted parametrized spectra or \emph{specters} for short.  Thus, a specter has the same relationship to a parametrized spectrum as a section of a line bundle has to a function.  Moreover, a specter for the haunt $\LL(c)$ associated to a class $c \in H^1(X;\PPic(\SS))$ is determined by a 0-cochain with coboundary $c$, that is by an element $f \in C^0(X;\OO_X)$ together with an identification $\delta f \cong c$.

This fundamental analogy is summarized in table~\ref{tableframework} and is explained in detail in the following sections.

\begin{sidewaystable}

\begin{tabular}{|c||c|c|}

\hline
&  \parbox[c][1.7cm]{7cm}{\begin{center} In Algebraic Geometry \end{center}}  &  
\parbox[c][1.7cm]{7cm}{\begin{center}In Twisted Parametrized \\ 
Stable Homotopy Theory\end{center}}
\\
\hline
\hline
The basic ring $R$ is  &  
\parbox[c][1.4cm]{6cm}{\begin{center} an ordinary ring \end{center}}  &  
\parbox[c][1.25cm]{6cm}{\begin{center} the category $\Sp$ of spectra \end{center}}
\\
\hline
The structure sheaf $\OO_X$ is given by  & 
\parbox[c][1.4cm]{6cm}{\begin{center} $\OO_X(U) = R$-valued \\
functions on $U$ \end{center}}  &  
\parbox[c][1.25cm]{6cm}{\begin{center} $\OO_X(U) =$ the category of \\
parametrized spectra on $U$ \end{center}} 
\\
\hline
The set of units in the ring $R$ is  &  
\parbox[c][1.4cm]{6cm}{\begin{center}the multiplicative group
$R^{\times}$ \end{center}}  &  
\parbox[c][1.25cm]{6cm}{\begin{center}the monoidal category $\Pic(\SS)$ \\ 
of invertible $\SS$-modules \end{center}} 
\\
\hline
The sheaf of units is given by &
\parbox[c][1.4cm]{6cm}{\begin{center} $\OO_X^{\times}(U) = R^{\times}$-valued \\
    functions on $U$ \end{center}} &
\parbox[c][1.25cm]{6cm}{\begin{center} $\PPic(\SS)(U) =$ the category
    of \\ invertible param. spectra on $U$ \end{center}} 
\\
\hline
A line bundle $\LL$ is  &  
\parbox[c][1.4cm]{6cm}{\begin{center} a sheaf that is an \\
invertible $\OO_X$-module \end{center}}  &  
\parbox[c][1.25cm]{6cm}{\begin{center} a stack that is a locally free \\ 
rank-one $\OO_X$-module \end{center}}
\\
\hline
Line bundles $\LL(c)$ are classified by  &  
\parbox[c][1.4cm]{6cm}{\begin{center} elements $c \in
    H^1(X;\OO_X^{\times})$ \end{center}} &
\parbox[c][1.25cm]{6cm}{\begin{center} elements $c \in
    H^1(X;\PPic(\SS))$ \end{center}}
\\
\hline
\parbox[c][1.25cm]{7cm}{\begin{center} A global section of a line bundle is \end{center}}  &  
\parbox[c][1.4cm]{6cm}{\begin{center} an ordinary section  \end{center}} &
\parbox[c][1.25cm]{6cm}{\begin{center} \emph{a twisted parametrized spectrum}
\end{center}}
\\
\hline
Such sections are classified by &
\parbox[c][1.4cm]{6cm}{\begin{center}elements $f \in C^0(X;\OO_X)$ \\
such that $\delta f = c$ \end{center}}  &  
\parbox[c][1.25cm]{6cm}{\begin{center}elements $f \in C^0(X;\OO_X)$ \\
together with $\delta f \cong c$ \end{center}}
\\
\hline
\end{tabular}
\vspace{10pt}
\caption{The analogy with basic algebraic geometry.} \label{tableframework}
\end{sidewaystable}

\subsection{The Structure Stack of Parametrized Spectra} \label{ssps}

We begin by describing the category of parametrized spectra and its associated homotopy theory.  A spectrum $E$ is, most naively, a series $E_i$ of based spaces equipped with structure maps $\Sigma E_i \ra E_{i+1}$ from the suspension of one space to the next.  Similarly, we can describe a parametrized spectrum over $X$ by giving a series $E_i$ of based spaces over $X$ together with structure maps $\Sigma_X E_i \ra E_{i+1}$ from the fibrewise suspension of one space to the next.  (A based space over $X$ is a space together with a projection map to $X$ and a section of this projection.)
This naive viewpoint is sufficient for many purposes, including taking the homology and cohomology of a base space $X$ with coefficients in a parametrized spectrum, but it fails to provide a foundation for a good smash product on the category of parametrized spectra.  As we are interested in considering this category to be a semi-ring, it is essential that we have a highly associative and commutative smash product.  We therefore work with the category of orthogonal parametrized spectra.  An orthogonal parametrized spectrum on $X$ is a diagram spectrum in the category of based spaces over $X$, where the diagram category is finite dimensional inner product spaces and their isometries.  See Mandell, et al.~\cite{mmss} for a description of diagram spectra and May-Sigurdsson~\cite{ms} for an extensive discussion of orthogonal parametrized spectra.  In keeping with our sheaf-theoretic philosophy, will we let $\OO_X(U)$ denote the category of (orthogonal) parametrized spectra on the open set $U \subset X$.

We have selected a model for parametrized spectra because we need to ensure that we have a well behaved smash product, but we are not of course interested in the point-set-level peculiarities of this particular model.  As such, we need to keep in mind a homotopy theory, that is a model structure or at least a notion of weak equivalences, on the category $\OO_X(U)$ of parametrized spectra; we should consider a category $\CC$ equivalent to $\OO_X(U)$ not if there is an equivalence of categories between them, but if there is a functor inducing an equivalence of homotopy theories.  We will focus on the homotopy theory on $\OO_X(U)$ associated to the stable model structure defined in~\cite{ms}.
Suffice it to say that the fibrant objects in the stable model structure are, in particular, quasi-fibrations of spectra over $X$, and the weak equivalences between fibrant objects are maps that induce weak equivalences on each fibre.  We can think of this as a homotopy theory of quasi-fibrations of spectra, rather than of all parametrized spectra; this reduction will be important in our consideration of local systems of categories of spectra in section~\ref{sectionlocsys}.  As an aside, we note that there may be other interesting model structures on parametrized spectra, for example ones in which there is a much larger class of fibrant objects; the formulation of twisted parametrized spectra in section~\ref{modstack} will work for these alternate model structures, producing a very different and perhaps even more intriguing theory.

\begin{remark}
We have fixed a notion of homotopy theory on the category of parametrized spectra, namely the one coming from the stable model structure.  The concept of an $\infty$-category conveniently encodes the notion of a category together with an associated homotopy theory.  An $\infty$-category is, roughly speaking, a category together with 2-morphisms, 3-morphisms, and so on, such that all the $n$-morphisms are invertible for $n>1$.  Though $\infty$-categories are as yet little utilized, many familiar structures, including simplicial categories, Segal categories, Segal spaces, and quasi-categories, give models for $\infty$-categories; see Lurie~\cite{lurietopoi,lurieinfty} for a thorough treatment.  The pairing of a category and a notion of homotopy theory will be so pervasive that in this section and in section~\ref{modstack} we will frequently use "category" to mean "$\infty$-category" and implicitly take associated notions, such as equivalence, monoidal structure, module, and so forth, to refer to their $\infty$-categorical analogs.  The reader who is bothered by the resulting inexplicitness should defer to section~\ref{sectionlocsys} where haunts and specters are recharacterized in more traditional terms.  
(Model-theoretically inclined readers may want to take "category" to mean "model category" and this will be perfectly suitable by way of understanding, but we will be utilizing categories of categories (read $\infty$-categories of $\infty$-categories) and the category of model categories is not known to have a model structure.)
\end{remark}

\twiddles

The association to an open set $U \subset X$ of the category $\OO_X(U)$ of parametrized spectra over $U$ is meant to function as a "sheaf of rings" analogous to the structure sheaf of $R$-valued functions on an algebraic variety.  First we consider the ring-like structure on $\OO_X(U)$ and then proceed to the sheaf- or stack-like properties of $\OO_X$.

The category $\OO_X(U)$ of parametrized spectra is a symmetric bimonoidal category in the sense of Laplaza~\cite{laplaza}; that is, it comes equipped with two symmetric monoidal functors (wedge and smash) and natural distributivity isomorphisms satisfying various coherence relations.  In fact, the category is better behaved that the average symmetric bimonoidal category because the wedge product is the categorical coproduct; the additive associativity isomorphisms and the distributivity isomorphisms are therefore canonically defined.  Of course, there is a rigidification functor~\cite{elmen-mand,dunn} that replaces a symmetric bimonoidal category with an equivalent bipermutative category (where the associativity isomorphisms are identity transformations).  We will frequently and implicitly use the bipermutative category associated to $\OO_X(U)$, particularly when discussing modules in the next section.  This symmetric bimonoidal or bipermutative structure makes $\OO_X(U)$, for all intents and purposes, into a semi-ring.

Philosophically, a stack $\CC$ is a presheaf of categories satisfying descent up to equivalence of categories; (see, for example, Moerdijk's treatment in~\cite{moerdijk}).  That is, the category $\CC(U)$ living over a large open set $U$ is determined, up to equivalence, by the categories $\CC(V)$ living over small open subsets $V \subset U$, in the same way that the value of a sheaf is determined by its local behavior.  Note that these categories $\CC(U)$ need not be groupoids (as is usually assumed) and that both "presheaf" and "equivalence" can be freely interpreted.  For example, "presheaf" might mean literal presheaf, presheaf up to coherent natural isomorphism, or presheaf up to coherent natural homotopy equivalence; similarly, "equivalence of categories" might mean ordinary equivalence of categories, Quillen equivalence of model categories, homotopical equivalence of $\infty$-categories, or something analogous.  In general, we take presheaf to mean presheaf up to coherent natural isomorphism, and equivalence of categories to mean homotopical equivalence of $\infty$-categories.

As the notation suggests, we have restriction functors $i^*:\OO_X(U) \ra \OO_X(V)$ associated to inclusions $i: V \subset U$.  Because these restrictions of parametrized spectra boil down to literal restriction functors in categories of topological spaces over the base $X$, these functors give $\OO_X$ the structure of a presheaf of categories on $X$.
This presheaf $\OO_X$ is a stack.  Though it is a stack in the usual, literal sense that it satisfies descent up to equivalence of categories (as can be checked using the prestack gluing condition given in~\cite[p.11]{moerdijk}), we are only concerned with the fact that it is a stack in the sense that it satisfies descent up to homotopical equivalence of categories.

\begin{summary}
The association $\OO_X$ to an open subset $U \subset X$ of the category $\OO_X(U)$ of orthogonal parametrized spectra on $U$ is a stack of symmetric bimonoidal ($\infty$-)categories.
\end{summary}

\subsection{Modules over the Structure Stack} \label{modstack}

We have replaced an ordinary ring $R$ by the semi-ring category $\Sp$ of spectra and we are investigating invertible sheaves in this new context.  We have introduced our basic "sheaf of rings" $\OO_X$, namely the structure stack of parametrized spectra.  In this section, we describe and classify locally free rank-one modules over this structure stack---we refer to these modules, briefly, as "haunts"---and we study their categories of global sections.  These global sections are the fundamental objects of twisted parametrized stable homotopy theory and we call them "twisted parametrized spectra" or "specters" for short.

\subsubsection{Haunts} \label{haunts}

A module over a stack $\RR$ of symmetric bimonoidal categories is a stack $\MM$ of symmetric monoidal categories together with an action $\RR \times \MM \ra \MM$ appropriately compatible (by analogy with a module over a ring) with the monoidal structures.  We do not spell out this compatibility; see Dunn~\cite{dunn} for an extensive discussion of modules over semi-ring categories, Lurie~\cite{lurieinfty} for some of the technicalities involved in semi-ring $\infty$-categories and their modules, and remark~\ref{coproductremark} below for an explanation of why we do not attend to the details of these compatibility relations. 
Such a module $\MM$ over a symmetric bimonoidal stack $\RR$ is locally free of rank one if for all points $x \in X$ there exists an open set $U \subset X$ containing $x$ and an object $S \in \MM(U)$ such that the map $\RR |_U \xra{} \MM |_U$ determined for $V \subset U$ by
\begin{align}
\RR(V) &\xra{} \MM(V) \nn\\
A &\mapsto A \cdot S |_V \nn
\end{align}
is an equivalence of symmetric monoidal stacks.
\begin{defn}
A \emph{haunt} on a space $X$ is a locally free rank-one module over the structure stack $\OO_X$ of parametrized spectra on $X$.
\end{defn}
\begin{remark} \label{coproductremark}
Because the additive monoidal structure in the category $\OO_X$ of parametrized spectra is given by the categorical coproduct, most of the compatibility conditions~\cite{dunn} for $\OO_X$-modules are automatically satisfied, provided the monoidal structure on the module is also the coproduct.  Indeed, it is generally sufficient to treat $\OO_X$ as a multiplicative monoid and study stacks with an action of this monoid.  As a point of philosophy, though, it is important to keep in mind that we are really dealing with modules over \emph{ring} stacks.
\end{remark}
\begin{remark}
We limit our attention to locally free rank-one modules over parametrized spectra, but we imagine that there may be quite interesting and intricate homotopy-theoretic information in the structure of higher rank modules.
\end{remark}

In doing geometry over an ordinary ring $R$, we think of invertible sheaves as line bundles.
Such a line bundle is most easily and explicitly described by taking trivial $R$-bundles over an open cover $\{U_{i}\}$ and specifying appropriate gluing data $r_{ij}: R \times U_{ij} \ra R \times U_{ij}$ on the two-fold intersections $U_{ij} = U_{i} \cap U_{j}$.  By analogy, we think of haunts as bundles whose fibre is the category $\Sp$ of spectra.
To specify such a bundle, we can take a trivial $\Sp$-bundle over an open cover and give gluing data $r_{ij}: \Sp \times U_{ij} \ra \Sp \times U_{ij}$ on the intersections; more precisely, this gluing data amounts to automorphisms $r_{ij}: \OO_X(U_{ij}) \ra \OO_X(U_{ij})$.  Before formalizing this viewpoint, we give two examples.

\begin{example} \label{s1haunt}
Cover the base space $X=\S^1$ by two open semicircles $U_0$ and $U_1$ and denote by $V$ and $W$ the two components of the intersection $U_{01}$.  Glue $\OO_X |_{U_0}$ and $\OO_X |_{U_1}$ together along $V$ by the identity map on $\OO_X(V)$ and along $W$ by the map
\begin{align}
\OO_X(W) &\ra \OO_X(W) \nn\\
T &\mapsto T \sm_W (\S^n \times W) \nn
\end{align}
In other words, the monodromy around the circle is the map $\Sp \ra \Sp$ given by suspension by $\S^n$.  Schematically the resulting haunt appears as in figure~\ref{circlehaunt}.  As we will see, every haunt over $\S^1$ is equivalent to this suspension haunt for some integer $n$.
\begin{figure}[h]
\begin{center}
\caption{A bundle over the circle with fibre the category of spectra} \label{circlehaunt}
\epsfig{figure=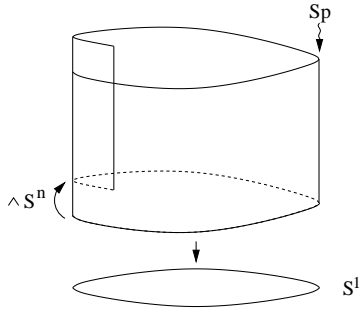,height=1.6in}
\end{center}
\end{figure}
\end{example}

\begin{example} \label{s3haunt}
Now take the base space $X$ to be $\S^3$ with its usual hemispherical cover by two open sets $D^3_0$ and $D^3_1$; the intersection of these open sets is the equatorial band $\S^2 \times (-\epsilon,\epsilon)$.  Let $\S^2 {\rtimes} \S^2$ denote the nontrivial $\S^2$ bundle over $\S^2$.  We can define a haunt over $X$ by gluing the trivial bundles $\Sp \times D^3_0$ and $\Sp \times D^3_1$ as follows:
\begin{align}
\OO_X(\S^2 \times (-\epsilon,\epsilon)) &\xra{} \OO_X(\S^2 \times (-\epsilon,\epsilon)) \nn\\
T &\mapsto T \sm_{\S^2 \times (-\epsilon,\epsilon)} ((\S^2 {\rtimes} \S^2) \times (-\epsilon, \epsilon)) \nn
\end{align}
This gluing produces the only nontrivial haunt over $\S^3$.
\end{example}

Line bundles are built by gluing together trivial $R$-bundles over intersections.  At a point, this gluing is determined by an $R$-module automorphism of the ring $R$, that is by an element of $\Aut_R R$.  Of course, any such automorphism is multiplication by a unit of $R$; in other words, $\Aut_R R \cong R^{\times}$.  The gluing data for a line bundle 
is therefore a 1-cocycle with values in $R^{\times}$ or more precisely, with values in the sheaf $\OO_X^{\times}$ of $R^{\times}$-valued functions.  Up to equivalence, these bundles are classified by the cohomology group $H^1(X;\OO_X^{\times})$.

Appropriately interpreted, all these facts remain true when $R$ is replaced by the category $\Sp$ of spectra---see table~\ref{tableframework}.  Any automorphism of $\Sp$ as a module over itself is given by smashing with an invertible spectrum; moreover, there is an equivalence of categories $\Aut_{\Sp} \Sp \cong \Sp^{\times}$.  Here the objects of $\Sp^{\times}$ are the invertible spectra and the morphisms are weak equivalences of spectra.  This category is denoted $\Pic(S^0)$ in the literature and we will use that notation to refer to both the category and its realization.  The gluing data for a haunt is a 1-cocycle with values in the sheaf $\PPic(S^0)$ of "$\Pic(S^0)$-valued functions".  This "sheaf" is a stack of monoidal categories; the category $\PPic(S^0)(U)$ has objects invertible parametrized spectra on $U$ and morphisms weak equivalences.  Naturally enough, haunts are classified up to equivalence by the cohomology group $H^1(X;\PPic(\SS))$.

The stack $\PPic(\SS)$ is locally constant in the sense that for any contractible $U$ the category $\PPic(\SS)(U)$ is weakly equivalent to the category $\Pic(\SS)$; as a result, the cohomology group $H^1(X;\PPic(\SS))$ is isomorphic to the group of homotopy classes of maps $[X,B\Pic(\SS)]$.  This isomorphism is a special feature of the homotopy-theoretic setting; in the discrete algebraic analog, the groups $H^1(X;\OO_X^{\times})$ and $[X,BR^{\times}]$ are quite distinct---the former classifies line bundles and the latter classifies line bundles with flat connection, that is local systems.  Indeed, there is no distinction in parametrized stable homotopy theory between a line bundle and a local system, and the local system perspective will form the framework for our discussion in section~\ref{sectionlocsys}.  In order to obviate the unpleasant details of  cohomology with coefficients in a stack of monoidal categories, we emphasize the classification of haunts in terms of homotopy classes of maps:

\begin{prop} \label{hauntclassification}
Haunts on $X$ are classified up to equivalence by the group of homotopy classes of maps from $X$ to $B\Pic(\SS)$.
\end{prop}

\nid This proposition will be established in detail in section~\ref{sectionlocsys}---see especially theorems~\ref{hauntclassp2} and~\ref{folkthm}.

Any invertible spectrum is weakly equivalent to a sphere of some integer dimension, so $B\Pic(S^0)$ has the homotopy type of $B(\Z \times B\GL_1(S^0))$ where $\GL_1(S^0)$ denotes the module automorphisms of the sphere spectrum, also known as the group of stable self equivalences of the sphere.  In particular, $\pi_1(B\Pic(S^0)) = \Z$ and $\pi_3(B\Pic(S^0)) = \Z/2$, explaining the classifications mentioned in examples~\ref{s1haunt} and~\ref{s3haunt}.  We will discuss the homotopy groups of $B\Pic(S^0)$ in more detail in the context of specter invariants in section~\ref{sectioninvariants}.

\subsubsection{Specters} \label{specters}

As previously mentioned, specters generalize parametrized spectra in the same way that sections of line bundles generalize functions:
\begin{defn}
A \emph{twisted parametrized spectrum} or \emph{specter} on $X$ is a global section of a haunt over $X$.  That is, it is an object $S \in \MM(X)$ of the category of global sections of a locally free rank-one module $\MM$ over the structure stack $\OO_X$ of parametrized spectra on $X$.
\end{defn}

An ordinary $R$-valued function on $X$ is determined by a 0-cocycle with values in the sheaf $\OO_X$ of $R$-valued functions, that is by a 0-cochain $f \in C^0(X;\OO_X)$ such that the coboundary vanishes: $\delta f = 0$.  Suppose $c \in Z^1(X;\OO_X^{\times})$ is a 1-cocycle defining a line bundle $L(c)$ with fibre $R$.  A section of $L(c)$ is presented by a 0-cochain $f \in C^0(X;\OO_X)$ cobounding the cocycle $c$, which is to say such that $\delta f = c$.   This section need not trivialize the line bundle because it is allowed to take non-invertible values, unlike the defining cocycle for the bundle.

Analogously, a parametrized spectrum on $X$ can be described by a 0-cocycle with values in the stack of parametrized spectra.  This amounts to giving a 0-cochain $f \in C^0(X;\OO_X)$, namely a parametrized spectrum on each open set of a cover, together with a compatible system of equivalences on intersections; this system of equivalences is concisely encoded in the equation $\delta f \cong 0$.
Let $c \in Z^1(X;\PPic(S^0))$ denote a 1-cocycle defining a haunt $L(c)$; concretely, this means that for a cover $\{U_{i}\}$ we have invertible parametrized spectra $c_{ij} \in \PPic(S^0)(U_{ij})$ on two-fold intersections together with fixed weak equivalences $\phi_{ijk}: c_{ij} \sm_{U_{ijk}} c_{jk} \ra c_{ik}$ satisfying the obvious coherence relation on four-fold intersections.  A specter for this haunt, that is a section of $L(c)$, is most easily presented by a 0-cochain $f \in C^0(X;\OO_X)$ together with an identification $\delta f \cong c$ of the coboundary of $f$ with the cocycle $c$.  What this means is that on two-fold intersections compatible equivalences are given between $f_i \sm_{U_{ij}} c_{ij}$ and $f_j$.  This cochain presentation is well suited to giving explicit examples of specters.

\begin{example} \label{s1specter}
Let $L_n$ denote the haunt over $S^1$ whose monodromy is suspension by $S^n$.  This haunt is depicted in figure~\ref{circlehaunt} and can be presented, roughly speaking, as follows: take two copies of the stack of parametrized spectra on an interval, that is of $\OO_{D^1}$, and glue them together at the two pairs of endpoints by the maps $\Sp \xra{\id} \Sp$ and $\Sp \xra{\sm S^n} \Sp$ respectively.  We now define a specter $T$ for $L_n$.  Over the first interval $D^1_0$ the specter is a trivial parametrized spectrum with fibre $S^n$; that is $T |_{D^1_0} = S^n \times D^1_0$.  Over the second interval $D^1_1$ the specter is a cone on $S^n \sqcup S^0$; that is, writing $D^1_1 = C(\pt \sqcup \pt)$, we have $T |_{D^1_1} = C(S^n \sqcup S^0)$.  See figure~\ref{circlespecter}.
\begin{figure}[h]
\begin{center}
\caption{A twisted parametrized spectrum over the circle} \label{circlespecter}
\epsfig{figure=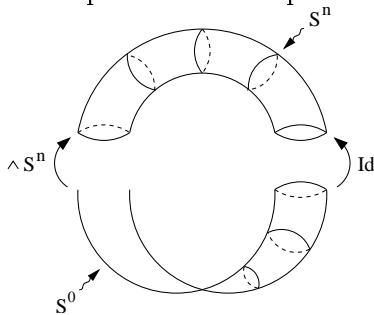,height=1.6in}
\end{center}
\end{figure}
\end{example}

\begin{example} \label{s3specter}
Let $L_{S^3}$ denote the nontrivial haunt over $S^3$ described in example~\ref{s3haunt}.  Roughly speaking, $L_{S^3}$ is constructed using the equatorial gluing function $\psi: \OO_{S^2} \ra \OO_{S^2}$
given by $\psi(P) = P \sm_{S^2} (S^2 {\rtimes} S^2)$.  Define a specter $T$ for $L_{S^3}$ by $T = (S^0 \times D^3) \cup_{\psi} C(S^2 {\rtimes} S^2)$; that is, on one hemisphere the specter is a trivial parametrized spectrum with fibre $S^0$ and on the other hemisphere (thought of as the cone $C(S^2)$) the specter is the cone on the bundle $S^2 {\rtimes} S^2$.
\end{example}

\begin{remark} \label{furutaremark} 
The procedure (illustrated in these examples) of piecing together parametrized spaces to give a global geometric object has also appeared in a preprint by Furuta~\cite{furuta} under the rubric "prespectra with parametrized universe".  Furuta's viewpoint is similar to ours in spirit, but his `parametrized prespectra' are substantially more rigid than specters; in particular, few specters can be realized as `parametrized prespectra'.  It would therefore seem to be difficult to develop a reasonable homotopy theory of twisted spectra (which is essential for applications to Floer homotopy) using the rigid geometry of `parametrized prespectra'.
\end{remark}

In section~\ref{ssps} we emphasized the fact that the homotopy theory of parametrized spectra should be conceived as a homotopy theory of fibrations of spectra.  The above examples of specters have singularities and are therefore not locally fibrations of spectra.  There is a fibrant replacement functor which takes such a singular specter and returns a locally fibrant specter with the same "global homotopy type"; (see section~\ref{sectioninvariants} for a discussion of global homotopy invariants of specters).  On the one hand, it is easier to explicitly describe and compute invariants of singular specters; on the other hand, it is easier to characterize and classify specters using fibrant presentations.  

The following model for fibrant replacement of specters bears close resemblance to the one May and Sigurdsson use for parametrized spectra; a more thorough technical treatment of the replacement functor can be found in their manuscript~\cite{ms}.  Suppose $S$ is a specter over the base $X$.  Let $PX$ denote the path fibration on $X$ and $s$ and $t$ the source and target maps $PX \ra X$.  A fibrant model $F(S)$ for the specter $S$ is very roughly given by $s_!(t^*(S))$---here $s_!$ denotes integration over the fibre in a sense analogous to that given for parametrized spectra in~\cite{ms}.  In other words, the points of the fibre $F(S)_x$ at a point $x \in X$ are pairs consisting of a path in $X$ from $x$ to $y$ and a point of the fibre $S_y$.  We proceed to some examples of fibrant specters.

\begin{example} \label{s1fibrant}
As before, let $L_n$ denote the haunt over $S^1$ whose monodromy is suspension by $S^n$.  Define a specter $T$ for $L_n$ as follows: the fibre $T_x$ of $T$ at every point $x \in S^1$ is $\bigvee_{i \in \Z} S^{n\cdot i}$ and the monodromy operator is the natural equivalence $\Sigma^n(\bigvee S^{n\cdot i}) \simeq \bigvee S^{n \cdot i}$.  This is a model for the specter in example~\ref{s1specter} which is locally a fibration of spectra.  More generally, any specter for the haunt $L_n$ on $S^1$ can be described by giving a spectrum $A$ together with an equivalence of $A$ with its n-th suspension:
\begin{equation} \nn
\{ \textrm{Specters} / (L_n, S^1) \} \leftrightsquigarrow \{ \textrm{Spectra } A \textrm{ with } \phi: \Sigma^n(A) \simeq A \}
\end{equation}
\end{example}

\begin{example} \label{s2fibrant}
There is only one nontrivial haunt over $S^2$; call it $L_{S^2}$.  Let $U$ and $V$ denote the two hemispheres of $S^2$.  The haunt $L_{S^2}$ is constructed, roughly speaking, by gluing $\OO_U$ and $\OO_V$ along the equatorial $S^1$ using the function $\psi:\OO_{S^1} \ra \OO_{S^1}$ given by $\psi(P) = P \sm_{S^1} (S^0 {\rtimes} S^1)$; here $S^0 {\rtimes} S^1$ denotes the nontrivial $S^0$ bundle over $S^1$.  Suppose we want to construct a specter $T$ for $L_{S^2}$ that is locally a fibration of spectra.  The restriction $T |_U$ of the specter to one hemisphere $U$ will be equivalent to the parametrized spectrum $A \times U$, for some spectrum $A$.  On the boundary of $U$, the specter is $A \times \partial U = A \times S^1$, and therefore on the boundary of $V$ the specter must be $\psi(A \times S^1) = (A \times S^1) \sm_{S^1} (S^0 {\rtimes} S^1)$---let us denote this last parametrized spectrum by $\A^{\tw(-1)}$.  As $\A^{\tw(-1)}$ is a fibrant parametrized spectrum over $S^1$ with fibre $A$, it is constructed by gluing together two copies of $A \times D^1$ along the boundaries $A \times S^0$.  Such a gluing is determined by a map $S^0 \ra \Aut(A)$, where $\Aut(A)$ denotes the homotopy automorphisms of the spectrum $A$.  In the case of $\A^{\tw(-1)}$ this gluing is the map $-1_A: S^0 \ra \Aut(A)$ taking one point to $\id_A$ and the other point to $-\id_A$.  The specter $T$ restricts on the hemisphere $V$ to a fibrant parametrized spectrum $T |_V$ with boundary $\A^{\tw(-1)}$---this parametrized spectrum defines a nullhomotopy of the gluing map $-1_A$ of $\A^{\tw(-1)}$.  In summary:
\begin{equation} \nn
\{ \textrm{Specters} / L_{S^2} \} \leftrightsquigarrow \{ \textrm{Spectra } A \textrm{ with } \phi: D^1 \ra \Aut(A) \textrm{ s.t. } \phi_0 = \id_A, \phi_1 = -\id_A \}
\end{equation}
\end{example}

\begin{example} \label{s3fibrant}
Recall the nontrivial haunt $L_{S^3}$ over $S^3$ determined by the equatorial gluing function $\psi(P)=P \sm_{S^2} (S^2 {\rtimes} S^2)$.  Note that this haunt is isomorphic to the haunt determined by the gluing function $\psi'(P)=P \sm_{S^2} (S^0 {\rtimes} S^2)$, where $S^0 {\rtimes} S^2$ is the unique nontrivial $S^0$-spectrum bundle over $S^2$.  Suppose $T$ is a specter for $L_{S^3}$ that is locally a fibration of spectra.  Then, as in the previous example, $T$ restricted to one hemisphere $U$ is a trivial bundle with fibre spectrum $A$.  The boundary $T |_{\partial U}$ of this restriction is $A \times \partial U = A \times S^2$ and the boundary $T |_{\partial V}$ of the other restriction must therefore be $\psi'(A \times S^2) = (A \times S^2) \sm_{S^2} (S^0 {\rtimes} S^2)$---denote this last parametrized spectrum by $\A^{\tw(\eta)}$.  The fibrant parametrized spectrum $\A^{\tw(\eta)}$ over $S^2$ is determined by the gluing function $\eta_A: S^1 \ra \Aut(A)$; here $\eta_A$ is the function $\eta \sm A$ where $\eta:S^1 \ra \Aut(S^0)$ is the nontrivial element of the first stable stem.  The restriction $T |_V$ of the specter to the second hemisphere constitutes a nullhomotopy of the gluing function $\eta_A$ for the boundary parametrized spectrum $\A^{\tw(\eta)}$.  Again we have a classification:
\begin{equation} \nn
\{ \textrm{Specters} / L_{S^3} \} \leftrightsquigarrow \{ \textrm{Spectra } A \textrm{ with } \phi:D^2 \ra \Aut(A) \textrm{ s.t. } \phi_{S^1} = \eta_A \}
\end{equation}
\end{example}

\twiddles

We have seen a variety of examples of specters and have described all the specters associated to a few particular haunts.  We now systematically investigate the equivalence classification of specters, giving a homotopy-theoretic description of specters that naturally parallels the characterization of haunts in proposition~\ref{hauntclassification}.

Line bundles are classified by $H^1(X;\OO_X^{\times})$, where $\OO_X^{\times}$ is the sheaf of invertible $R$-valued functions; in section~\ref{haunts}, we noted that haunts are classified by $H^1(X;\PPic(\SS))$, where $\PPic(\SS)$ is the stack of invertible parametrized spectra.  Specters for the haunt $L(c)$ associated to a class $c \in H^1(X;\PPic(\SS))$ are classified by 0-cochains $f \in C^0(X;\OO_X)$ together with an identification $\delta f \cong c$.  We remarked that the group $H^1(X;\PPic(\SS))$ is isomorphic to $[X,B\Pic(\SS)]$ and that the latter is a simpler characterization of the classifying group for haunts; as before $\Pic(\SS)$ denotes the realization of the category of invertible spectra.  We describe a classification of specters based on the $[X,B\Pic(\SS)]$ classification of haunts.

We take our cue, as usual, from the discrete algebraic analogue.  The space $BR^{\times}$ is a classifying space for ordinary line bundles \emph{with flat connection}.  There is a universal $R^{\times}$-bundle $ER^{\times}$ over $BR^{\times}$ and an associated $R$-bundle $P(R):=ER^{\times} \times_{R^{\times}} R$.  A flat section of the line bundle $L(c)$ associated to the map $c: X \ra BR^{\times}$ is determined by a lift of $c$ to a map $s: X \ra P(R)$:
\begin{equation} \nn
\xymatrix{
R \ar[r] & P(R) \ar[d] \\
X \ar[r]_-{c} \ar@{-->}[ur]^-{s} & BR^{\times}
}
\end{equation}

This description of sections of a line bundle with flat connection easily translates into the context of haunts and specters.  The space $B\Pic(\SS)$ is a classifying space for haunts.  Let $\Sp_w$ denote the realization of the subcategory of weak equivalences of the category of spectra, and note that $\Pic(\SS)$ acts on $\Sp_w$.  There is a universal $\Pic(\SS)$-bundle $E\Pic(\SS)$ over $B\Pic(\SS)$ and an associated $\Sp_w$-bundle \mbox{$P(\Sp_w):=E\Pic(S^0) \times_{\Pic(S^0)} \Sp_w$}.  We can consider lifts of a classifying map $c: X \ra B\Pic(\SS)$ to $P(\Sp_w)$:
\begin{equation} \nn
\xymatrix{
{\Sp}_w \ar[r] & P({\Sp}_w) \ar[d] \\
X \ar[r]_-{c} \ar@{-->}[ur] & B\Pic(S^0)
}
\end{equation}
Indeed, such lifts classify specters:
\begin{prop} \label{specterclassification}
Let $c:X \ra B\Pic(S^0)$ be the classifying map for a haunt $L(c)$.  Weak equivalence classes of specters for the haunt $L(c)$ are in one-to-one correspondence with homotopy classes of lifts of $c$ to maps $X \ra P(\Sp_w)$.
\end{prop}
\nid That there should be such a homotopy-theoretic classification of specters was suggested to us by Bill Dwyer.  A slightly stronger result will be stated and proved as theorem~\ref{spectclass-simp}; granting theorem~\ref{folkthm}, this proposition follows directly from corollary~\ref{specterclass-space}.

We conclude our initial discussion of haunts and specters with a few remarks about products.  Given two ordinary invertible sheaves (line bundles) $L$ and $L'$ over $X$ we can form their tensor product $L \otimes_{\OO_X} L'$.  Analogously, given two haunts we should be able to form their tensor product.  Making sense of tensoring two modules over a stack of semi-ring categories would require a bit of doing; (see Dunn~\cite{dunn} for a definition of the tensor product of modules over a bipermutative category).  
Using the classification of haunts in proposition~\ref{hauntclassification} we can side step this categorical tensor construction: define the product $L \otimes L'$ of two haunts to be the haunt classified by the product (in the group structure on $B\Pic(S^0)$) of the classifying maps $c,c':X \ra B\Pic(S^0)$.  
There is a product of specters covering this tensor product of haunts.  The classifying projection $P(\Sp_w) \ra B\Pic(S^0)$ for specters is a map of multiplicative monoids.  Given two specters $T$ and $T'$ classified by lifts $s,s':X \ra P(\Sp_w)$ of the haunt maps $c,c':X \ra B\Pic(S^0)$, we simply define the product $T \sm T'$ to be the specter (for the haunt $L \otimes L'$) classified by the product lift $s \cdot s':X \ra P(\Sp_w)$.

\section{Local Systems of Categories of Spectra}   \label{sectionlocsys} 

In section~\ref{sectionsheaves}, we described haunts as invertible sheaves, or more specifically as locally free rank-one modules over the structure stack of parametrized spectra; specters, the twisted parametrized spectra, were global sections of these modules.  In this section, we leave behind that sheaf-theoretic approach and reformulate haunts as local systems of categories of spectra.

An ordinary local system on a variety $X$ over the ring $R$ is a line bundle with flat connection.  The flat sections of this bundle form a sheaf $L$ with the property that the ring of sections $L(U)$ is isomorphic to $R$ for any sufficiently small open set $U$.  This sheaf encodes all the data of the local system, and we can therefore present any local system as follows: to each open set in a cover $\{U_i\}$ assign a ring $C_i$ isomorphic to $R$, and to the two-fold intersections assign a family of compatible isomorphisms $f_{ij}: C_i \xra{\cong} C_j$.  Any automorphism of $R$ is given by multiplication by an element of $R^{\times}$; local systems are therefore classified by the first cohomology group $H^1(X;R^{\times})$.  Note that this group is quite distinct from the group $H^1(X;\OO_X^{\times})$ classifying invertible sheaves---this difference boils down to the fact that no matter how small the open set $U$, the ring $\OO_X(U)$ of $R$-valued functions on $U$ is not isomorphic to $R$ as a ring.

Our favorite "ring" is, of course, the category $\Sp$ of spectra, and the presentation of local systems of categories of spectra is analogous to that of ordinary local systems: to each open set in a cover $\{U_i\}$ assign a category $C_i$ weakly equivalent to the category of spectra, and to the
two-fold intersections assign a family of compatible equivalences $f_{ij}: C_i \xra{\simeq} C_j$.  A self-equivalence of the category of spectra is given by smashing with an element of $\Pic(\SS)$, the category of invertible spectra; haunts, now conceived of as local systems of categories of spectra, are therefore classified by the group $H^1(X;\Pic(\SS))$.  This group \emph{is} isomorphic to the group classifying invertible sheaves of categories of spectra and indeed the notions of invertible sheaf and local system are interchangeable in the context of twisted parametrized stable homotopy theory.  Underlying this correspondence is the fact that the category $\OO_X(U)$ of parametrized spectra on a contractible open set $U$ is Quillen equivalent to the category $\Sp$ of nonparametrized spectra.

An invertible sheaf is a line bundle and a local system is a line bundle with flat connection; associated to a local system there is therefore an underlying invertible sheaf.  The sections of a local system are the locally constant or \emph{flat} sections of the associated invertible sheaf.  Analogously, a local system of categories of spectra has an underlying invertible sheaf of categories of spectra.  The sections of the local system are in particular sections of the corresponding invertible sheaf; indeed, these sections are locally constant in the sense that they are locally quasifibrations of spectra, and so are precisely the \emph{fibrant} twisted parametrized spectra. 

The fundamental analogy between local systems in algebraic geometry and in twisted parametrized stable homotopy theory is summarized in table~\ref{secondtableframework} and elaborated throughout the following sections.

\begin{sidewaystable}

\begin{tabular}{|c||c|c|}

\hline
&  \parbox[c][1.7cm]{7cm}{\begin{center} In Algebraic Geometry \end{center}}  &
\parbox[c][1.7cm]{7cm}{\begin{center}In Twisted Parametrized \\
Stable Homotopy Theory\end{center}}
\\
\hline
\hline
The basic ring $R$ is  &
\parbox[c][1.7cm]{6cm}{\begin{center} an ordinary ring \end{center}}  &
\parbox[c][1.25cm]{6cm}{\begin{center} the category $\Sp$ of spectra \end{center}}
\\
\hline
\parbox[c][1.7cm]{7cm}{\begin{center} A line bundle with \\ flat
    connection is \end{center}} &
\parbox[c][1.25cm]{6cm}{\begin{center} a local system \end{center}} &
\parbox[c][1.25cm]{6cm}{\begin{center} a diagram of categories, \\
each weakly equivalent to $\Sp$ \end{center}}
\\
\hline
\parbox[c][1.7cm]{7cm}{\begin{center} Line bundles with flat \\
    connection are classified by \end{center}} &
\parbox[c][1.25cm]{6cm}{\begin{center} elements $c \in H^1(X;
    R^{\times})$ \\ or maps $X \ra BR^{\times}$ \end{center}} &
\parbox[c][1.25cm]{6cm}{\begin{center} elements $c \in
    H^1(X;\Pic(\SS))$ \\ or maps $X \ra B\Pic(\SS)$ \end{center}}
\\
\hline
\parbox[c][1.7cm]{7cm}{\begin{center} The global sections of a line
    bundle \\ with flat connection are \end{center}} &
\parbox[c][1.25cm]{6cm}{\begin{center} the flat sections
\end{center}} &
\parbox[c][1.25cm]{6cm}{\begin{center} the fibrant twisted \\ parametrized
    spectra \end{center}}
\\
\hline
\parbox[c][1.7cm]{7cm}{\begin{center} Such sections are classified by \end{center}} &
\parbox[c][2.5cm]{6cm}{\begin{center} lifts
\parbox[c][0cm]{3.8cm}{\begin{center}
$
\xymatrix{ & ER^{\times} \times_{R^{\times}} R \ar[d] \\
X \ar[r] \ar@{-->}[ur] & BR^{\times}
}
$
\end{center}}
\end{center}} &
\parbox[c][2.5cm]{6cm}{\begin{center} lifts
\parbox[c][0cm]{5.1cm}{\begin{center}
$
\xymatrix{ & E\Pic(\SS) \times_{\Pic(\SS)} \Sp \ar[d] \\
X \ar[r] \ar@{-->}[ur] & B\Pic(\SS)
}
$
\end{center}}
\end{center}}
\\
\hline
\end{tabular}
\vspace{10pt}
\caption{The analogy with basic algebraic geometry, continued.} \label{secondtableframework}
\end{sidewaystable}

\twiddles

As described above, a local system of categories of spectra is determined by associating to each open set in a cover a category weakly equivalent to the category of spectra and to intersections of open sets appropriate equivalences of categories.  In order to formalize this viewpoint, we need a "category of categories of spectra"---such a category is typically referred to as a homotopy theory of homotopy theories.  In section~\ref{htpyhtpy} we discuss a convenient choice of such a category, namely the model category of simplicial categories.

In section~\ref{diagsimp} we use simplicial categories to make precise the notion of a local system of categories of spectra.  Suppose the base space $X$ is the realization of a simplicial complex $B$, and the cover of $X$ is the cover of $B$ by the stars of its simplicies.  Then a haunt on $X$, thought of as a local system, is a functor from the category of simplicies of $B$ to the category of simplicial categories, such that the functor takes objects to categories weakly equivalent to the category of spectra and such that the functor takes morphisms to weak equivalences.  This simplicial category definition of a haunt has the disadvantage that it clouds the conceptual simplicity and the geometry of haunts and their associated specters, but it has the advantage that it avoids the technicalities of stacks of $\infty$-categories and thereby eases the proofs of propositions~\ref{hauntclassification} and~\ref{specterclassification}---these propositions follow immediately from the results in sections~\ref{haunts2},~\ref{specters2}, and~\ref{automspec}.

In section~\ref{sectionthom}, we give \emph{another} reformulation of haunts and specters: haunts are $\AA_{\infty}$ ring spectra arising as Thom spectra of multiplicative stable spherical fibrations on loop spaces, and the associated specters are simply modules over the ring spectrum.  This last description has the advantage of being entirely elementary---it avoids both modules over stacks of parametrized spectra and diagrams in the model category of simplicial categories---but it thoroughly obscures various constructions with and applications of specters, and so in section~\ref{sectioninvariants} we return to our original sheaf-theoretic perspective.

\subsection{The Homotopy Theory of Homotopy Theories} \label{htpyhtpy}

As it will play a central role in our discussion of local systems, we briefly describe the homotopy theory of homotopy theories.  In recent decades, model categories have been the predominant notion of abstract homotopy theory.  However, there is not known to be a model structure on the category of model categories, and this is a huge impediment to constructing (as we are doing in this paper) bundles of homotopy theories.  We must consider a weaker notion of abstract homotopy theory in order to have a decent homotopy theory of homotopy theories.  There are various options, including the model category of Segal categories (due to Hirschowitz-Simpson~\cite{hir-simp}), the model category of complete Segal spaces (due to Rezk~\cite{rezk}), and the model category of simplicial categories (due to Dwyer-Hirschhorn-Kan~\cite{dhk} and Bergner~\cite{bergner}).  Which we pick does not matter because all three model categories are Quillen equivalent (a result due to Bergner~\cite{bergner}); we work with the model category of simplicial categories, as this is the simplest to describe.

By a `simplicial category' we will mean a category enriched over simplicial sets.  Dwyer and Kan~\cite{dk80b} realized that to a model category $M$ there is canonically associated a simplicial category $L^H M$, the hammock localization of $M$, which encodes all of the homotopy-theoretic information contained in $M$.  This is the sense in which simplicial categories are a faithful representation of abstract homotopy theories.  In particular, the homotopy category $\Ho(M)$ of $M$ is recovered as the category of components $\pi_0(L^H M)$ of the hammock localization; here the category of components $\pi_0(C)$ of a simplicial category $C$ (also called the homotopy category of $C$) has the same objects as $C$, but has morphisms $\Hom_{\pi_0(C)}(x,y) = \pi_0(\Hom_C(x,y))$.  A morphism $f:x \ra y$ in a simplicial category $C$ is called a homotopy equivalence if it becomes an isomorphism in $\pi_0(C)$.

The model structure on the category of simplicial categories is as follows.  A map $\phi: C \ra D$ of simplicial categories is a weak equivalence if it is a Dwyer-Kan equivalence, namely if $\phi$ is a weak equivalence on $\Hom$ sets and an equivalence on homotopy categories; that is, $\phi$ is a weak equivalence if $\phi: \Hom_C(x,y) \ra \Hom_D(\phi(x),\phi(y))$ is a weak equivalence of simplicial sets for all objects $x,y\in C$, and if $\phi: \pi_0(C) \ra \pi_0(D)$ is an equivalence of categories.  A map $\phi: C \ra D$ is a fibration if it is a fibration on $\Hom$ sets and if all homotopy equivalences in $D$ lift to $C$; that is, $\phi$ is a fibration if $\phi:\Hom_C(x,y) \ra \Hom_D(\phi(x),\phi(y))$ is a fibration of simplicial sets for all $x,y \in C$, and if for all objects  $x \in C$ and all homotopy equivalences $h:\phi(x) \ra z$ in $D$, there exists a homotopy equivalence $\widetilde{h}:x \ra y$ in $C$ such that $\phi(\widetilde{h})=h$.
Cofibrations of simplicial categories are determined, as usual, by the left lifting property.  These define a model structure on the category of simplicial categories~\cite{bergner}.  Given an object in a model category, there is a good notion of the space of automorphisms of that object, and in the following we will be focused on the automorphisms of the category of spectra (thought of as a simplicial category via its hammock localization).

\begin{remark} 
Considering that our objection to model categories as a representation of homotopy theory was that there is no obvious model category of model categories, it seems odd to have insisted on having a model category of simplicial categories rather than merely a simplicial category of simplicial categories.  Of course, we can recover a simplicial category of simplicial categories as the hammock localization of the model category of simplicial categories, but it would be better not to have to rely on the crutch of a model structure.  The real solution to this and many other problems is to work directly in the $\infty$-category of $\infty$-categories.  We do not do this because the details of such an $\infty$-categorical theory are not yet fully in place; see Lurie~\cite{lurietopoi,lurieinfty}, though, for substantial progress in that direction.
\end{remark}

\subsection{Diagrams of Simplicial Categories} \label{diagsimp}

We described haunts as locally free rank-one modules over the structure stack of parametrized spectra.  We now reinterpret this notion in terms of diagrams of simplicial categories weakly equivalent to the hammock localization of the category of spectra.  We then express the category of specters for such a haunt as a homotopy limit in the model category of simplicial categories.

\subsubsection{Haunts} \label{haunts2}

Let $H$ be a locally free rank-one module over the structure stack on $X$.  By definition, if $U$ is a sufficiently small open set in $X$, then the category $H(U)$ is homotopy equivalent to the category $\OO_X(U)$ of parametrized spectra on $U$.  If the subspace $U$ is moreover contractible, then $\OO_X(U)$ is homotopy (indeed Quillen) equivalent to the category $\Sp$ of spectra.  (Roughly speaking, this is because a fibrant parametrized spectrum is a quasifibration and any quasifibration over a contractible space is trivializable---see~\cite{ms}.)  Let $\{U_i\}$ denote a contractible cover of $X$; the module $H$ determines, by the above remarks, an assignment to each $U_i$ a category $C_i$ equivalent to the category of spectra and to the two-fold intersections $\{U_{ij}\}$ a collection of compatible equivalences $C_i \simeq C_j$---this is the data of a local system of categories of spectra.  We formalize the compatibilities in the collection of equivalences using diagram functors into the category of simplicial categories.

First we fix some notation.  The space $X$ is homotopy equivalent to the realization of a simplicial set $B$; let $s(B)$ denote the category of simplicies of $B$, that is the category whose objects are the simplicies of $B$ and whose morphisms are the face and degeneracy maps.  Let $\sCat$ denote the model category of simplicial categories, and let $\Sp$ now denote the object of $\sCat$ given by the hammock localization of the category of spectra.  Let $w(\sCat,\Sp)$ denote the weak equivalence component of $\sCat$ containing $\Sp$; that is, the objects of $w(\sCat,\Sp)$ are simplicial categories that can be connected to $\Sp$ by a zig-zag of weak equivalences, and the morphisms are weak equivalences of simplicial categories.

\begin{remark} 
 We pause to consider a few set-theoretic issues.  The category $\sCat$ is really the category of small simplicial categories.  The hammock localization of the category of spectra is not only not small, it need not even have small $\Hom$ sets; we therefore chose a small simplicial category homotopically equivalent to that hammock localization---by the notation $\Sp$ we will implicitly refer to that small replacement.  In a similar vein, the weak equivalence component of $\Sp$ in $\sCat$ is not a small category; we will need to use it in constructions that only apply to small categories, so we chose a small subcategory of this weak equivalence component that is homotopically equivalent to the full component---we implicitly refer to that small replacement by the notation $w(\sCat,\Sp)$.  We will not henceforth distinguish between such (simplicial) categories and their small replacements.
\end{remark} 

\begin{rechar}
A \emph{haunt} over a simplicial set $B$ is a functor from $s(B)$, the category of simplicies of $B$, to $w(\sCat,\Sp)$, the weak equivalence component of the category of simplicial categories containing the category of spectra.  The \emph{category of haunts} over $B$, denoted $\Haunt_B$, is the full diagram category $w(\sCat,\Sp)^{s(B)}$.
\end{rechar} 

We immediately have a notion of the \emph{space of haunts}, namely the realization $|N.\Haunt_B|$ of the nerve of this diagram category.  Note that we will not in general distinguish between simplicial sets and their realizations.  There is a natural candidate for a classifying space for haunts, namely $N.w(\sCat,\Sp)$.  The idea that a weak equivalence component of an object of a model category can function as a classifying complex is of course due to Dwyer and Kan~\cite{dk84}.  Indeed, there is a suggestive homotopy equivalence $N.w(\sCat,\Sp) \simeq B\haut(\Sp)$.  Here $B\haut(\Sp)$ is the nerve of the simplicial category with one object and with morphisms the simplicial monoid $\haut(\Sp)$ of homotopy automorphisms of the category of spectra; (this simplicial monoid is defined to be the sub-simplicial monoid of $\Hom_{L^H(\sCat)}(\Sp,\Sp)$ consisting of the components projecting to isomorphisms in $\Hom_{\pi_0(L^H(\sCat))}(\Sp,\Sp)$).  The classification of haunts can therefore be expressed as follows.

\begin{theorem} \label{hauntclassp2}
The space $N.\Haunt_B$ of haunts over a simplicial set $B$ is weakly homotopy equivalent to the (derived) mapping space $\Hom(B,B\haut(\Sp))$.  In other words, $B\haut(\Sp)$ is a classifying space for haunts.
\end{theorem}

\begin{remark}
We adopt the convention that all mapping spaces are implicitly derived unless otherwise noted.  We will also take "holim" and "hocolim" to refer to the homotopically invariant homotopy limit and colimit functors; these are sometimes referred to as the corrected homotopy limit and colimit and can be defined respectively by composing functorial objectwise fibrant or cofibrant replacement with the Bousfield-Kan holim or hocolim functor.  
\end{remark}

\nid Before proving the theorem, we state one lemma:

\begin{lemma}
Suppose $M$ is a model category that is Quillen equivalent to a cofibrantly generated simplicial model category.  Let $B$ be a simplicial set and let $s(B)$ denote the category of simplicies of $B$.  For any object $X$ of $M$ there is a weak homotopy equivalence
\begin{equation} \nn
N.(w(M,X)^{s(B)}) \simeq \holim\limits_{s(B)} N.w(M,X)
\end{equation}
\end{lemma}

\nid Dwyer and Kan prove this for $M$ equal to the category of simplicial sets~\cite[Thm 3.4]{dk84}, but
their proof works for any cofibrantly generated simplicial model category.  Moreover, both sides of the equivalence are weakly homotopy invariant under Quillen equivalence between not-necessarily-simplicial model categories; this follows using various results from~\cite{dk80a,dk80b}.

\begin{proof}[Proof of Theorem~\ref{hauntclassp2}]
We have the chain of equivalences:
\begin{align}
N.\Haunt_B \equiv N.(w(\sCat,\Sp)^{s(B)}) &\simeq \holim\limits_{s(B)} N.w(\sCat,\Sp) \nn\\
&\simeq \Hom(N.s(B),N.w(\sCat,\Sp)) \nn\\
&\simeq \Hom(B,B\haut(\Sp)) \nn
\end{align}
The category of simplicial categories is Quillen equivalent to the category of complete Segal spaces~\cite{bergner} which is a cofibrantly generated simplicial model category~\cite{rezk}; the first equivalence therefore follows from the above lemma.  The second equivalence is a consequence of~\cite[Prop 18.2.6]{hirsch}, and the third is immediate.
\end{proof}

\subsubsection{Specters} \label{specters2}

We now discuss specters in this new context of diagrams of simplicial categories.  In the sheaf-theoretic framework of section~\ref{sectionsheaves} a specter was a global section of a haunt.  A haunt is now a functor $H$ from the category of simplicies $s(B)$ into the category $\sCat$ of simplicial categories.  Naturally enough, a "section" of such a diagram $H$ of simplicial categories should be some appropriately consistent choice of objects $\{x_{b} \in H(b)\}_{b \in s(B)}$ of the simplicial categories $H(b)$ in the diagram---we can think of an object $x_{b} \in H(b)$ as a locally constant (or flat) section of the haunt $H$ restricted to the simplex $b$.  It would be too much to ask that the collection of objects $x_{b}$ be strictly compatible with the morphisms in the diagram $H$.  Instead, we merely demand that there be chosen homotopies to "glue the objects together"---this gluing data is formally encoded in a homotopy limit.

\begin{rechar}
Let $H$ be a haunt over $B$, that is a functor from the category of simplicies $s(B)$ to the category of simplicial categories that lands in the weak equivalence component $w(\sCat,\Sp)$ of the category of spectra.  The (simplicial) \emph{category of specters} for the haunt $H$, denoted $\Specter_H$, is defined to be the homotopy limit $\holim\limits_{s(B)} H$.
\end{rechar}

Of course, this homotopy limit is only defined because we have a model structure on the category of simplicial categories, and as usual we mean the homotopically invariant homotopy limit.  We also have an associated \emph{space of specters} for the haunt $H$, namely $N.w(\Specter_H)$.  Here $w(C)$ denotes the sub-simplicial category of the simplicial category $C$ whose objects are the same as those of $C$ but whose morphisms $\Hom_{w(C)}(a,b)$ are the components of $\Hom_C(a,b)$ projecting to isomorphisms in $\Hom_{\pi_0(C)}(a,b)$; note that $N.w(C)$ is, a priori, a bisimplicial set and we implicitly take its diagonal.

There isn't a classifying space for specters, per se, but there is a classifying fibration; that is, there is a fibration $\psi: U{\haut(\Sp)} \ra B\haut(\Sp)$ such that the space of specters for a fixed haunt $H:B \ra B\haut(\Sp)$ is homotopy equivalent to the space of lifts of $H$ along $\psi$.  The fibre of $\psi$ should be the "space of spectra" and we think of the total space $U{\haut(\Sp)}$ as the "universal haunt".   In fact, this classifying fibration comes from a fibration of simplicial categories, which is defined as the diagonal map in the following diagram:

\begin{equation} \nn
\xymatrix{
\hocolim\limits_{w(\sCat,\Sp)} (-) \ar[d] 
\ar@{^{(}->}[r]^{\sim} 
& 
\overline{\hocolim\limits_{w(\sCat,\Sp)} (-)}  
\ar@{->>}[dl]^{\Psi} 
\\
\hocolim\limits_{w(\sCat,\Sp)} \pt &
}
\end{equation}

\nid That is, we factor the left hand map by a weak equivalence followed by a fibration $\Psi$; this fibration is the desired specter classifying fibration of simplicial categories.  Here $\hocolim\limits_{w(\sCat,\Sp)} (-)$ refers to the homotopy colimit of the inclusion $w(\sCat,\Sp) \ra \sCat$.  

The analogous specter classifying fibration of spaces is the right hand vertical arrow in the diagram

\begin{equation} \nn
\xymatrix{
N.w\left(\overline{\hocolim\limits_{w(\sCat,\Sp)} (-)}\right) \ar@{^{(}->}[rr]^-{\sim} \ar@{->>}[d] && U\haut(\Sp) \ar@{->>}[d]^{\psi} \\  
N.w\left(\hocolim\limits_{w(\sCat,\Sp)} \pt\right) \ar@{<->}[r]^-{\simeq} & N.w(\sCat,\Sp) \ar@{<->}[r]^-{\simeq} & B\haut(\Sp)
}
\end{equation}

\nid The fibration factorization here defines the space $U\haut(\Sp)$.  Note that the fibre of $\psi$ over the point $X \in w(\sCat,\Sp)$ is weakly equivalent to $N.w(X)$ which is in turn weakly equivalent to the space of spectra $N.w(\Sp)$.  
Morally speaking, the bundle $U{\haut(\Sp)}$ is the bundle $E\haut(\Sp) \times_{\haut(\Sp)} N.w(\Sp)$ associated to the tautological bundle $E\haut(\Sp) \ra B\haut(\Sp)$.  However, on its face the action of the simplicial monoid $\haut(\Sp)$ on $N.w(\Sp)$ is only defined up to weak homotopy, which is insufficient for defining the associated bundle; presumably the action can be made strict, but we do not pursue that here.

We now state the simplicial-category-level classification result for specters:

\begin{theorem} \label{spectclass-simp}
Let $H: s(B) \ra w(\sCat,\Sp)$ be a functor defining a fixed haunt on the simplicial set $B$.  This functor determines an associated map $h:\hocolim\limits_{s(B)} \pt \ra \hocolim\limits_{w(\sCat,\Sp)} \pt$ classifying the haunt.  The category of specters $\Specter_H$ for this haunt is weakly equivalent, as a simplicial category, to the category of lifts of $h$ along the specter classifying fibration:
\begin{equation} \nn
\xymatrix{
& \overline{\hocolim\limits_{w(\sCat,\Sp)} (-)} \ar@{->>}[d]^{\Psi} 
\\
\hocolim\limits_{s(B)} \pt \ar@{-->}[ur] \ar[r]_{h} & \hocolim\limits_{w(\sCat,\Sp)} \pt
}
\end{equation}
In other words the category of specters is weakly equivalent to a (derived) mapping space in the overcategory of $\hocolim\limits_{w(\sCat,\Sp)} \pt$, namely
\begin{equation} \nn
\Specter_H \simeq \Hom_{\hocolim\limits_{w(\sCat,\Sp)} \pt}\left(\hocolim\limits_{s(B)} \pt, \hocolim\limits_{w(\sCat,\Sp)} (-)\right)
\end{equation}
\end{theorem}

\begin{proof}
The first step in the proof is a consequence of the following general lemma:
\begin{lemma}
Let $M$ be either an $\infty$-topos (such as simplicial sets or spaces) or a model category of homotopy theories (such as simplicial categories, complete Segal spaces, or quasi-categories).  Let $D$ be a small category and $F:D \ra M$ a functor.  Denote by $\phi: \hocolim\limits_D F \ra \hocolim\limits_D \pt$ the natural projection.  Provided $F$ takes all morphisms in $D$ to weak equivalences in $M$, the homotopy limit of $F$ is weakly equivalent to the object of derived sections of the map $\phi$; that is
\begin{equation} \nn
\holim\limits_D F \simeq \Hom_{\hocolim\limits_D \pt}\left(\hocolim\limits_D \pt, \hocolim\limits_D F\right)
\end{equation}
\end{lemma}
\nid In the case of specters, we therefore have the equivalence
\begin{equation} \nn
\Specter_H \equiv \holim\limits_{s(B)} H \simeq \Hom_{\hocolim\limits_{s(B)} \pt}\left(\hocolim\limits_{s(B)} \pt,
\hocolim\limits_{s(B)} H\right)
\end{equation}
We want to translate this mapping space into a space of lifts of the classifying map $h: \hocolim\limits_{s(B)} \pt \ra \hocolim\limits_{w(\sCat,\Sp)} \pt$.  We begin by rewriting one of the homotopy colimits in terms of a larger indexing category:
\begin{equation} \nn
\hocolim\limits_{s(B)} H \simeq h^*\left(\hocolim\limits_{w(\sCat,\Sp)} (-)\right)
\end{equation}
Here $h^*$ denotes the derived pullback from the overcategory in $\sCat$ of $\hocolim\limits_{w(\sCat,\Sp)} \pt$ to the overcategory of $\hocolim\limits_{s(B)} \pt$.  Next, we have a Quillen adjunction
\begin{equation} \nn
\xymatrix{
\sCat / \left(\hocolim\limits_{s(B)} \pt\right) \ar@<-.5ex>[r]_{h_!} & \sCat / \left(\hocolim\limits_{w(\sCat,\Sp)} \pt\right) \ar@<-.5ex>[l]_{h^*}
}
\end{equation}
where the pushforward $h_!$ is given by precomposition with the map $h$.  This adjunction leads to an equivalence of function complexes
\begin{equation} \nn
\Hom_{\hocolim\limits_{s(B)} \pt}\left(\hocolim\limits_{s(B)} \pt, h^*\left(\hocolim\limits_{w(\sCat,\Sp)} (-)\right)\right) \simeq
\Hom_{\hocolim\limits_{w(\sCat,\Sp)} \pt}\left(\hocolim\limits_{s(B)} \pt, \hocolim\limits_{w(\sCat,\Sp)} (-)\right)
\end{equation}
as desired.
\end{proof}

Not surprisingly, the analogous result at the level of spaces is the following.

\begin{cor} \label{specterclass-space}
Let $h: B \ra B\haut(\Sp)$ denote the classifying map for a fixed haunt $H$ over the simplicial set $B$.  The associated space of specters $N.w\Specter_H$ is weakly homotopy equivalent to the space $\Hom_{B\haut(\Sp)}(B,U\haut(\Sp))$ of maps from $B$ to the universal haunt $U\haut(\Sp)$ that commute with projection to $B\haut(\Sp)$.  In other words, the space of specters for the haunt $H$ is weakly equivalent to the space of lifts in the diagram
\begin{equation} \nn
\xymatrix{
& U{\haut(\Sp)} \ar[d]^{\psi} & N.w(\Sp) \ar[l] \\
B \ar@{-->}[ur] \ar[r]_-{h} & B\haut(\Sp) &
}
\end{equation}
\end{cor}

\begin{proof}
The chain of equivalences is
\begin{align}
N.w\Specter_H = N.w\holim\limits_{s(B)} H &\simeq N.w \Hom_{\hocolim\limits_{w(\sCat,\Sp)} \pt}\left(\hocolim\limits_{s(B)} \pt, \hocolim\limits_{w(\sCat,\Sp)} (-)\right) \nn\\
& = N.w \Hom_{\hocolim\limits_{w(\sCat,\Sp)} \pt}\left(\hocolim\limits_{s(B)} \pt, \overline{\hocolim\limits_{w(\sCat,\Sp)} (-)}\right) \nn\\
& \simeq N.\Hom_{\hocolim\limits_{w(\sCat,\Sp)} \pt}\left(\hocolim\limits_{s(B)} \pt, w\left(\overline{\hocolim\limits_{w(\sCat,\Sp)} (-)}\right)\right) \nn\\
& \simeq \Hom_{N.w(\sCat,\Sp)} \left(N.s(B), N.w\left(\overline{\hocolim\limits_{w(\sCat,\Sp)} (-)}\right)\right) \nn\\
& \simeq \Hom_{B\haut(\Sp)}(B,U\haut(\Sp)) \nn
\end{align}.

\nid The first homotopy equivalence is a consequence of the theorem, and the second line follows from the definition of the $\Hom$ set as a derived mapping space.  Next note that the functor $w:\sCat \ra \hGpd$ from simplicial categories to homotopy groupoids is right adjoint to the inclusion---the homotopy equivalence in the third line follows because $\hocolim\limits_{s(B)} \pt$ and $\hocolim\limits_{w(\sCat,\Sp)} \pt$ are already homotopy groupoids.  The category of homotopy groupoids is in fact Quillen equivalent to the category of simplicial sets; the equivalence in the fourth line follows, and the fifth line is immediate.
\end{proof}

\subsection{$\AA_{\infty}$ Thom Spectra on Loop Spaces} \label{sectionthom}

We begin this section by identifying the monoid $\haut(\Sp)$ of automorphisms of the category of spectra with a classifying space $\Z \times BG$ for stable spherical fibrations.  Using this identification, we can associate to a haunt over $X$ an $\AA_{\infty}$ ring spectrum arising as a Thom spectrum over the loop space of $X$.  This in turn allows us to recharacterize the category of specters for the haunt as a category of modules over that ring spectrum.  This recharacterization on the one hand obscures the intrinsic symmetry of specters and breaks the natural connection with parametrized homotopy theory, which is essential to the definition of specter invariants in section~\ref{sectioninvariants}; on the other hand, because module spectra are familiar objects, the change in perspective demystifies specters and will be important in applications to symplectic Floer homotopy~\cite{douglashopkins}.

\subsubsection{Automorphisms of the Category of Spectra} \label{automspec}

Haunts over a space $X$ are classified by maps from $X$ to the space $B\haut(\Sp)$, a deloop of the simplicial monoid of homotopy automorphisms of the category of spectra.  We investigate the homotopy type of this classifying space.  The category of spectra has a natural monoidal structure, the smash product; given an invertible spectrum $J$, the functor $J \sm -: \Sp \ra \Sp$ smashing with $J$ determines a self homotopy equivalence of the category of spectra.  Roughly speaking, this association determines a map from the category of invertible spectra (which we called $\Pic(S^0)$ in section~\ref{sectionsheaves}) to the space of self equivalences $\haut(\Sp)$.  That this map is a weak equivalence is well known to experts, but we are not aware of a statement or a proof in the literature:

\begin{theorem} \label{folkthm}
Let $\Pic(S^0)$, the Picard category, denote the subcategory of the category of spectra whose objects are invertible spectra and whose morphisms are weak equivalences.  There is a weak equivalence
\begin{equation} \nn
\Pic(S^0) \simeq \haut(\Sp)
\end{equation}
from the nerve of the Picard category to the simplicial set of self homotopy equivalences of the category $\Sp$ of spectra.
\end{theorem}

\begin{proof} We merely sketch the proof.  The category of spectra is a model category representing a particular homotopy theory, and we can work in any of a number of equivalent categories of homotopy theories.  We have primarily utilized $\sCat$, the category of simplicial categories, but for this theorem it is more convenient to work in $\qCat$, the category of quasi-categories---there is a Quillen equivalence between $\sCat$ and $\qCat$~\cite{lurieinfty}.  Recall that a quasi-category is a simplicial set that satisfies a weak Kan condition, namely that a horn $\partial \Delta^n \backslash \Delta^{n-1}_i$ fills in provided the missing face $\Delta^{n-1}_i$ is internal, that is $0<i<n$; this weak Kan condition reflects the idea that morphisms (edges) in a category are composable, but need not be invertible up to homotopy.  By $\qCat$ we refer, in fact, to the category of simplicial sets equipped with a model structure in which quasicategories are precisely the fibrant objects.

Let $\Sp$ denote a quasicategory modeling the category of spectra; we presume that $\Sp$ is equipped with a monoidal structure modeling the smash product.  In this context $\Pic(S^0)$ is a subquasicategory of $\Sp$ which is described as follows.  The vertices of $\Pic(S^0)$ are the invertible objects in $\Sp$, that is the vertices $v \in \Sp$ such that there exists a $w \in \Sp$ with $v \sm w$ weakly equivalent to $S^0 \in \Sp$; the $k$-simplicies of $\Pic(S^0)$ are the $k$-simplicies of $\Sp$ all of whose vertices are invertible and all of whose edges are weak equivalences.  By definition, the simplicial monoid $\haut(\Sp)$ has $k$-simplicies the set of weak equivalences $\Delta^k \times \Sp \xra{\simeq} \Sp$.
There is now a natural map
\begin{equation} \nn
\mu: \Pic(S^0) \ra \haut(\Sp)
\end{equation}
which takes a $k$-simplex $P$ in $\Pic(S^0)_k$ to the composite $\Delta^k \times \Sp \xra{P \times \id} \Pic(S^0) \times \Sp \xra{} \Sp \times \Sp \xra{} \Sp$.  This composite is an equivalence and is therefore a $k$-simplex in $\haut(\Sp)$.

We would like $\mu$ to be an equivalence.  It suffices to show that any map $F:(\Delta^k,\partial\Delta^k) \ra (\haut(\Sp),\Pic(S^0))$ is homotopic, relative to its boundary, to a map $F':(\Delta^k,\partial\Delta^k) \ra (\Pic(S^0),\Pic(S^0))$.  Suppose $k=0$, so $F$ is simply an equivalence $\Sp \ra \Sp$; we take $F'$ to be the map $\Sp \ra \Sp$ given by smashing with $F(S^0)$, that is, $F' = F(S^0) \in \Pic(S^0)$.  One extremely convenient feature of quasicategories (as distinguished from, for example, model categories) is that one can naturally take homotopy colimits over any simplicial set, not only over a category;  this feature is helpful in defining a comparison map between $F$ and $F'$.  The quasicategory $\Sp$ is in particular a simplicial set, and given an object $X \in \Sp$, let $(-/X)$ denote the "overcategory" of $X$, that is the subsimplicial set of $\Sp$ whose $0$-simplicies are maps $Y \ra X$.  Moreover, denote by $(-/X)_{\textrm{sph}}$ the corresponding "spherical subcategory", that is the full subcategory whose $0$-simplicies are the maps $S^i \ra X$.  We have the comparison map
\begin{equation} \nn
F(S^0) \sm X \simeq \hocolim\limits_{(-/X)_{\textrm{sph}}} F(S^i) \ra \hocolim\limits_{(-/X)} F(Y) \simeq F(X)
\end{equation}
This map is an equivalence when $X = S^0$, and the left hand side preserves homotopy colimits.  It is a consequence of Lurie's extensive work on quasicategories~\cite{lurieinfty} that because $F$ is an equivalence, it preserves homotopy colimits.  The comparison map is therefore a natural weak equivalence, as desired.  The cases of higher $k$ could be handled similarly.
\end{proof}

\begin{cor}
The simplicial set $\haut(\Sp)$ of self equivalences of the category of spectra has the homotopy type $\Z \times BG$ where $G$ is the space of stable self homotopy equivalences of the sphere, that is $G = \haut(\SS) \simeq \colim\limits_n \haut(S^n)$ where in the last expression $S^n$ denotes the ordinary n-sphere.
\end{cor}

\begin{proof}
Given the theorem, this is a consequence of the weak equivalence $\Pic(S^0) \simeq \Z \times BG$.  To see that equivalence, first note that any invertible spectrum is weakly equivalent to some shift $S^n$ of the sphere spectrum.  Thus the category $\Pic(S^0)$ has $\Z$ components; the n-th component is all spectra weakly equivalent to $S^n$ together with all weak equivalences between them.  By Dwyer and Kan's classification theorem~\cite{dk84}, this component has the homotopy type $B\haut(S^n) \simeq B\haut(\SS)$.
\end{proof}  

\subsubsection{Specters as Module Spectra} \label{specmod}

Armed with the identification of the classifying space $B\haut(\Sp)$ with $B(\Z \times BG)$, we can describe the $\AA_{\infty}$ ring spectrum corresponding to a haunt.  Let $h:X \ra  B\haut(\Sp) \simeq B(\Z \times BG)$ be the classifying map for a haunt over the space $X$.  The map $\Omega h: \Omega X \ra \Z \times BG$ defines a stable spherical fibration, which we will denote $\eta(h)$, over $\Omega X$.  
Because $\Omega h$ is a loop map, the spherical fibration $\eta(h)$ is multiplicative and the associated Thom spectrum $\Th(\eta(h))$ is therefore an $\AA_{\infty}$ ring spectrum~\cite{mahowald}.  The fibration $\eta(h)$ can be thought of more geometrically as follows.  The haunt is a local system or bundle over $X$ whose fibre is the category of spectra, and the monodromy of the haunt around a loop $\ell \in \Omega X$ is an invertible spectrum, namely the fibre $\eta(h)_{\ell}$; the multiplicative structure of the spherical fibration corresponds, naturally enough, to the composition of the loop monodromies.

If the Thom spectrum $\Th(\eta(h))$ encodes the structure of the haunt $h$, we might expect to be able to describe the associated category of specters in terms of this Thom spectrum.  Let $\TT$ be a specter for the haunt $h$ and suppose $\TT$ is fibrant in the sense that it is locally isomorphic to a quasifibration of spectra.  Given a loop $\ell:S^1 \ra X$ in $X$, we can pull $\TT$ back to a specter $\ell^*\TT$ on $S^1$. From example~\ref{s1fibrant} we know that this specter is determined by giving the spectrum $T$ at the basepoint of $S^1$ together with an equivalence $\phi_{\ell}: \eta(h)_{\ell} \sm T \simeq T$; in other words, we need to glue $T$ back to itself, but shifted by the monodromy sphere $\eta(h)_{\ell}$ along the given loop.  The family of compatible equivalences $\{\phi_{\ell}\}_{\ell \in \Omega X}$ amounts precisely to an action of  $\Th(\eta(h))$ on $T$.  To a specter for the haunt $h$ we can therefore associate a module over the ring spectrum $\Th(\eta(h))$; indeed there is an equivalence of categories:

\begin{prop}
Let $H$ denote a fixed haunt over $X$ with classifying map $h:X \ra B(\Z \times BG)$.  The loop map $\Omega h$ determines a multiplicative stable spherical fibration $\eta(h)$ with associated Thom spectrum $\Th(\eta(h))$ an $\AA_{\infty}$ ring spectrum.  There is a weak equivalence of simplicial categories
\begin{equation} \nn
\Specter_H \simeq L^H(\Th(\eta(h))\textrm{-mod})
\end{equation}
between the category of specters for $H$ and the hammock localization of the category of module spectra over the Thom spectrum $\Th(\eta(h))$.
\end{prop}

\nid That specters can be thought of as modules over a ring spectrum was also realized by Mike Hopkins and Jeff Smith, and the formulation here owes various details to discussions with them.

\begin{proof}[Sketch of proof]
We have already done most of the work in establishing, in theorem~\ref{spectclass-simp}, that
\begin{equation} \nn
\Specter_H \simeq \Hom_{\hocolim\limits_{w(\sCat,\Sp)} \pt}\left(\hocolim\limits_{s(B)} \pt, \hocolim\limits_{w(\sCat,\Sp)} (-)\right)
\end{equation}
Here $B$ is a simplicial set with the homotopy type of $X$.  By theorem~\ref{folkthm} the base space $\hocolim\limits_{w(\sCat,\Sp)} \pt \simeq N.w(\sCat,\Sp)$ has the homotopy type $B\Pic(\SS)$.  In particular, we can model this space by the simplicial category, also denoted $B\Pic(\SS)$, with just one object $\pt$ and with morphism space the (simplicial) monoid $\Pic(\SS)$.  Similarly, the total space $\hocolim\limits_{w(\sCat,\Sp)} (-)$ of the classifying fibration can be explicitly modeled by a simplicial category, denoted $U\Pic(\SS)$, as follows.  The objects of $U\Pic(\SS)$ are ordinary cofibrant and fibrant spectra.  The morphism space $\Hom_{U\Pic(\SS)}(T,S)$ of $U\Pic(\SS)$ is the homotopy colimit of the functor $\Hom_{\Sp}(- \sm T, S): \Pic(\SS)_{cf} \ra \sSet$.  Here $\Pic(\SS)_{cf}$ denotes the category of cofibrant and fibrant invertible spectra.  The idea behind this construction is that roughly speaking a morphism from $T$ to $S$ in $U\Pic(\SS)$ should consist of a pair $(\gamma,\phi)$ of an invertible spectrum $\gamma \in \Pic(\SS)$ and a morphism of spectra $\phi: \gamma \sm T \ra S$.  
Because $\hocolim\limits_{\Pic(\SS)_{cf}} \pt \simeq \Pic(\SS)_{cf}$, the projection $\Hom_{\Sp}(- \sm T, S) \ra \pt$ induces a map $U\Pic(\SS) \ra B\Pic(\SS)$.  Thinking of $X \simeq \hocolim\limits_{s(B)} \pt$ as a simplicial category, we therefore have the reformulation:
\begin{equation} \nn
\Hom_{\hocolim\limits_{w(\sCat,\Sp)} \pt}\left(\hocolim\limits_{s(B)} \pt, \hocolim\limits_{w(\sCat,\Sp)} (-)\right)
\simeq \Hom_{B\Pic(\SS)}(X,U\Pic(\SS))
\end{equation}
The final equivalence is
\begin{equation} \nn
\Hom_{B\Pic(\SS)}(X,U\Pic(\SS)) \simeq L^H(\Th(\eta(h))\textrm{-mod})
\end{equation}
On this count we merely indicate the map from the right to the left hand side.  We can model $X$ by the simplicial category with one object and with morphism space $\Omega X$.  The map $h: X \ra B\Pic(\SS)$ corresponds to the map of morphism spaces $\Omega X \ra \Pic(\SS)$ classifying the spherical fibration $\eta(h)$.  Given a $\Th(\eta(h))$-module $T$, one immediately has, for any $k$-simplex $s:\Delta^k \ra \Omega X$, a morphism $\Th(\eta(h) |_s) \sm T \ra T$.  By a shift in perspective, the stable spherical bundle $\eta(h)  |_s$ on $\Delta^k$ naturally corresponds to a $k$-simplex in the morphism category $\Pic(\SS)$ of $B\Pic(\SS)$, and the map $\Th(\eta(h) |_s) \sm T \ra T$ then provides a lift of this morphism $k$-simplex to $U\Pic(\SS)$, as desired.
\end{proof}

\begin{example} \label{s1module}
Consider again the haunt $L_n$ over $S^1$ whose monodromy is suspension by $S^n$.  This haunt is classified by the map $S^1 \ra B(\Z \times BG)$ representing $n \in \Z = \pi_1(B(\Z \times BG))$.  The loop of this map, $\Omega S^1 \simeq \Z \xra{n} \Z \ra \Z \times BG$, classifies the fibration over $\Z$ whose fibre at $i$ is $S^{n \cdot i}$.  The associated Thom spectrum is $\bigvee_{i \in \Z} S^{n \cdot i}$, which we considered as a specter already in example~\ref{s1fibrant}.  The set of specters for $L_n$ is, as we saw in that example, the same as the set of module spectra over this distinguished specter $\bigvee_{i \in \Z} S^{n \cdot i}$.
\end{example}

\section{Invariants of Specters}  \label{sectioninvariants}  

A parametrized spectrum $P$ over a space $X$ has two naturally associated spectra, namely the total spectrum $P/X$ representing the homotopy type of $P$ and the spectrum of sections $\Gamma(P)$ representing the cohomotopy type of $P$.  The generalized homology groups of these associated spectra provide invariants of the parametrized spectrum.  A specter, that is a twisted parametrized spectrum, has no globally defined homotopy or cohomotopy type analogous to the total spectrum $P/X$ or the spectrum of sections $\Gamma(P)$.  Nevertheless, it is frequently possible to define invariants associated to a specter by first applying a generalized homology functor and \emph{then} taking an associated global spectrum, rather than vice versa as is typical in parametrized homotopy theory.

We return to the stack-theoretic perspective of section~\ref{sectionsheaves}.  There, we treated the category of spectra as a ring and defined haunts to be locally free rank-one modules over a parametrized version of this ring; specters were global sections of these modules.  We begin this section by describing an analogous construction where the basic ring is the category of $R$-modules for a commutative ring spectrum $R$.  This leads to a notion of $R$-haunt and $R$-specter.  We then see how to associate to a haunt $H$ an $R$-haunt $H_R$ and to a specter $T$ for $H$ an $R$-specter $T_R$ for $H_R$; this base change is the aforementioned "generalized homology functor".  When the $R$-haunt $H_R$ is trivializable, the $R$-specter $T_R$ has the form of a parametrized $R$-module and therefore has a global homotopy type $T_R / X$.  The homotopy groups of $T_R / X$ are the $R$-homology invariants of the original specter $T$.  We describe a few examples of these specter invariants and discuss a spectral sequence for computing them.

\subsection{$R$-Haunts and $R$-Specters} \label{Rhaunt}

The category of spectra, or equivalently the category of modules over the sphere spectrum $\SS$, has a smash product which, roughly speaking, gives it the structure of a commutative ring.  The stack of parametrized spectra on a space $X$ therefore functions as a sheaf of rings, and we characterized haunts as locally free rank-one modules over this stack.  There are natural subcategories of the category of spectra that have their own commutative products, and we study modules over the stacks associated to these categorical rings.  Specifically, if $R$ is an $\AA_{\infty}$ ring spectrum, then we have the category $\Rmod$ of modules over $R$.  If $R$ is moreover $\E_{\infty}$, then there is a natural product $\sm_R$ which, roughly speaking, gives $\Rmod$ the structure of a commutative ring; see Dunn~\cite{dunn} for a detailed discussion of the monoidal structures on such module categories.  We
can associate to an open set $U \subset X$ the category $\OO^R_X(U)$ of parametrized $R$-modules on $U$.  These parametrized $R$-modules form a stack, which again has a monoidal structure coming from $\sm_R$.  Morally, an $R$-haunt is a locally free rank-one module over this structure stack $\OO^R_X$ of parametrized $R$-modules.  The stack $\OO^R_X$ is naturally a stack of $\infty$-categories, and by a module over this stack we refer to a stack of $\infty$-categories with an appropriate action of $\OO^R_X$.  A more complete definition of $R$-haunt is therefore as follows:

\begin{defn}
An \emph{$R$-haunt} is a stack $\MM$ of $\infty$-categories on a space $X$ with an action of the monoidal stack of $\infty$-categories $\OO_X^R$ of parametrized $R$-modules satisfying the following condition: the stack $\MM$ is locally free of rank one in the sense that for all points $x \in X$ there exists an open set $U \subset X$ containing $x$ and an object $Q \in \MM(U)$ such that the map $\OO_X^R  |_U \xra{} \MM  |_U$ given by $A \mapsto A \cdot Q$ is an equivalence of stacks of $\infty$-categories.
A twisted parametrized $R$-module or \emph{$R$-specter} is a global section of an $R$-haunt.
\end{defn}

In order to make this precise, one must give a thorough treatment of stacks of $\infty$-categories and of monoidal stacks of $\infty$-categories; we do not do this, but note that Lurie's work~\cite{lurieinfty} provides key elements of such a treatment.  Also note that, as in remark~\ref{coproductremark}, the additive structure on these stacks is given by the categorical coproduct and so is justifiably ignored---this saves us the horror of contemplating symmetric bimonoidal stacks of $\infty$-categories.

An $R$-haunt is locally equivalent to the stack of parametrized $R$-modules and as such can be specified concretely in terms of gluing functions.  Suppose $\{U_i\}$ is an open cover of $X$; on an intersection $U_{ij}$ a gluing function is a self equivalence of $\OO^R_X  |_{U_{ij}}$ as a module over itself.  Such an equivalence is given by smashing with an invertible parametrized $R$-module.  Let $\Pic(R)$ denote the category of invertible $R$-modules together with their homotopy equivalences, and let $\PPic(R)$ denote the corresponding sheaf of invertible parametrized $R$-modules.  The gluing data for an $R$-haunt is therefore a 1-cocycle $c$ with values in $\PPic(R)$, which is to say a compatible system $c_{ij}$ of parametrized $R$-module gluing functions on the one-fold intersections of the cover.  As we might guess from the classification of ordinary haunts in proposition~\ref{hauntclassification}, $R$-haunts on $X$ are classified by homotopy classes of maps from $X$ to $B\Pic(R)$.

An $R$-specter looks locally like a parametrized $R$-module, but has a global twist determined by the $R$-haunt.  Such a twisted parametrized $R$-module can be presented as follows: to describe an $R$-specter for the $R$-haunt associated to a gluing function $c$ for the cover $\{U_i\}$ it suffices to give a parametrized $R$-module $f_i$ on each open set $U_i$ together with compatible equivalences between $f_i \sm_{(R,U_{ij})} c_{ij}$ and $f_j$.  Compare section~\ref{specters}.  There is a classification of $R$-specters analogous to that of ordinary specters in proposition~\ref{specterclassification}, but we do not go into detail.

\twiddles

In section~\ref{haunts2} we formulated haunts not as stacks of $\infty$-categories, but as diagrams in the category of $\infty$-categories---or more precisely in the category $\sCat$ of simplicial categories.  The diagram had the homotopy type of the base space $X$, and the functor to $\sCat$ took objects to simplicial categories weakly equivalent to the category of spectra, and morphisms to weak equivalences.  A similar approach is possible for $R$-haunts but it requires more work.  Specifically, let $\Rmod$ denote the $\infty$-category of $R$-modules; this is an element of the category $\infCat$ of $\infty$-categories.  Now define $(\Rmod)\textrm{-mod}$ to be the subcategory of $\infCat$ of modules over $\Rmod$; in turn $w((\Rmod)\textrm{-mod},\Rmod)$ denotes the subcategory of $(\Rmod)\textrm{-mod}$ of objects weakly equivalent to $\Rmod$, together with the $(\Rmod)$-module weak equivalences between them.  Finally, an $R$-haunt would be a functor from an appropriate diagram homotopy equivalent to $X$ into $w((\Rmod)\textrm{-mod},\Rmod)$.  The associated $\infty$-category of $R$-specters would be the homotopy limit of this functor.  

There is also a Thom spectrum approach to $R$-haunts and $R$-specters.  Suppose an $R$-haunt $H_R$ is classified by the map $h_R:X \ra B\Pic(R)$.  The loop of this map $\Omega h_R: \Omega X \ra \Pic(R)$ classifies a fibration $\eta(h_R)$ of invertible $R$-modules over $\Omega X$; (indeed, the homotopy type of $\Pic(R)$ is $\Pic^0(R) \times B\GL_1(R)$ where $\Pic^0(R)$ denotes the equivalence classes of invertible $R$-modules, and $\GL_1(R)$ denotes the $R$-module self-equivalences of $R$).  This fibration $\eta(h_R)$ is in particular a parametrized spectrum and has an associated total spectrum $\eta(h_R) / X$ which we denote suggestively $\Th(\eta(h_R))$.  This total spectrum is an associative $R$-algebra and it encodes the structure of the $R$-haunt.  An associated $R$-specter is simply a $\Th(\eta(h_R))$-module, by which we mean an $R$-module $M$ together with an appropriately compatible action $\Th(\eta(h_R)) \sm_R M \ra M$.

\subsection{$R$-Homology of Specters} \label{Rhomology}

We now describe how to associate to haunts and specters respectively $R$-haunts and $R$-specters and we discuss the basic construction of specter invariants.  There is a natural map $\Smod \ra \Rmod$ from the category of spectra to the category of $R$-modules, given by $T \mapsto T \sm_{\SS} R$.  This is a map of rings and it underlies the fundamental base-change operation from specters to $R$-specters.  Given a haunt, that is a locally free rank-one module $\MM$ over the stack of parametrized spectra $\OO_X$, we can form the tensor stack $\MM \otimes_{\OO_X} \OO_X^R$ where $\OO_X^R$ is the stack of parametrized $R$-modules.  This tensor will be an $R$-haunt.  (See Dunn~\cite{dunn} for a definition of tensoring over ring categories.)  Any global section $P \in \MM(X)$ of $\MM$ transforms to the section $P \otimes R$; we therefore have an $R$-specter associated to any ordinary specter.  

In more down-to-earth terms, a haunt is presented by gluing together trivial bundles using invertible parametrized $\SS$-modules $c_{ij}$; the gluing functions for the associated $R$-haunt are simply $c_{ij} \sm R$.  Similarly, a specter is locally given by parametrized spectra $f_i$ and the associated $R$-specter is presented by the parametrized $R$-modules $f_i \sm R$.  Implicit here is the fact that the map $\Smod \ra \Rmod$ restricts to a map $\Pic(\SS) \ra \Pic(R)$ of Picard groups and this latter map deloops to a map $B\Pic(\SS) \ra B\Pic(R)$ of classifying spaces.  Thus we also have the purely homotopy-theoretic characterization of the haunt transformation, namely that the $R$-haunt associated to the haunt $X \ra B\Pic(\SS)$ is classified by the composite $X \ra B\Pic(\SS) \ra B\Pic(R)$.

\begin{remark}
The map $\Smod \ra \Rmod$ given by smashing with $R$ makes sense for any $\AA_{\infty}$ ring spectrum $R$.  This map induces a map $\GL_1(\SS) \ra \GL_1(R)$ and even a map $B\GL_1(\SS) \ra B\GL_1(R)$.  However, in order to build the $R$-haunt associated to an ordinary haunt, we need moreover a map $B(\Z \times B\GL_1(\SS)) \ra B(\Pic^0(R) \times B\GL_1(R))$.  Even barring the issue of what $\Pic^0(R)$ should mean, the space $B\GL_1(R)$ cannot deloop unless $R$ is commutative.  This provides another indication that in order to have $R$-homology invariants of specters, $R$ must be a commutative ring spectrum.  This commutativity requirement on the generalized homology invariants might appear surprising and like a fluke of the formulation, but in fact it reflects an essential aspect of the mathematical structure; later on we will see that the slogan is `semi-infinite homotopy types only have commutative generalized homology invariants'.
\end{remark}

Given a specter $P$ for the haunt $H$, we have "taken its $R$-homology" and produced an $R$-specter $P_R$ for the $R$-haunt $H_R$.   The essential idea behind specter invariants is that the structure of this $R$-haunt $H_R$ might be substantially simpler than that of $H$.  In particular, if $H_R$ is trivializable, that is if the composite $X \ra B\Pic(\SS) \ra B\Pic(R)$ is null homotopic, then any trivialization $\tau: H_R \xra{\sim} (X \times \Rmod)$ transforms the $R$-specter $P_R$ into a parametrized $R$-module $\tau(P_R)$.  The associated total spectrum $\tau(P_R) / X$ represents a global $R$-homotopy type for the specter $P$, even though $P$ does not itself have a global homotopy type.  The homotopy groups of $\tau(P_R) / X$ are what we might call the $R$-homology groups of $P$, denoted $R_*^{\tau}(P)$; these are the most straightforward and most easily computable invariants of specters.  Note that these groups definitely do depend on the trivialization $\tau$, but this ambiguity can be identified and controlled.

\begin{summary}
Let $P$ be a specter for the haunt $H$ and suppose that $\tau$ is a trivialization of the associated $R$-haunt $H_R$.  Then the \emph{$R$-homology groups of the specter $P$} are defined to be the homotopy groups of the total spectrum of the trivialization of the associated $R$-specter $P_R$:
\begin{equation} \nn
R_i^{\tau}(P) := \pi_i(\tau(P_R) / X)
\end{equation}
\end{summary}

The potential ambiguity in the trivialization of a trivializable $R$-haunt $H_R$ is governed by the space of automorphisms of the trivial $R$-haunt $X \times \Rmod$.  More specifically, the space of trivializations of $H_R$ is a torsor for $\Aut(X \times \Rmod) \simeq \Hom(X,\Pic(R))$.  Homotopic trivializations $\tau$ and $\tau'$ determine the same invariants $R^{\tau}_*(P) \simeq R^{\tau'}_*(P)$, and so we need only consider the set of components of the space of trivializations of $H_R$---this set of components is a torsor for $[X,\Pic(R)]=[X,\Pic^0(R) \times B\GL_1 R]$.  It often happens that $[X,B\GL_1 R]$ has only one element; if $X$ is connected, the set of components of the space of trivializations is then a torsor for $\Pic^0(R)$.  In this situation, which is to say when we have a specter $P$ for a haunt $H$ with $H_R$ trivializable and $[X,\Pic(R)]=\Pic^0(R)$, we can describe how the $R$-homology of $P$ is affected by a change in trivialization as follows.  A given trivialization $\tau$ can be modified by an invertible $R$-module $M \in \Pic^0(R)$ to the trivialization $M \cdot \tau$, and we have 
\begin{equation} \nn
R_*^{M \cdot \tau}(P) = \pi_*((M \cdot \tau)(P_R)/X) = \pi_*((\tau(P_R) \sm_R M)/X) = \pi_*(\tau(P_R)/X \sm_R M).
\end{equation}
If the ring spectrum $R$ is such that the homotopy $\pi_*(M)$ of an invertible module $M$ is projective over $R_*$, then we conclude that $R_*^{M \cdot \tau}(P) = R_*^{\tau}(P) \otimes_{R_*} \pi_*(M)$---in other words, changing the trivialization shifts the homology of $P$ by the homotopy of an invertible $R$-module.  This shift ambiguity is always present; we think of a homology group determined up to such a shift as uniquely determined and we ignore the $\tau$-dependency, writing simply $R_*(P)$.

\twiddles

To have homology invariants of a specter $P$ over a haunt $H$, we need to find a ring spectrum $R$ such that the composite $X \ra B\Pic(\SS) \ra B\Pic(R)$ is nullhomotopic.  As such, we need a thorough understanding of the homotopy types of various $B\Pic(R)$ and of the transformations $B\Pic(\SS) \ra B\Pic(R)$.  We already noted that $B\Pic(R) \simeq B(\Pic^0(R) \times B\GL_1(R))$.  By definition, $\GL_1(R)$ consists of the unit components of the zero space of the spectrum $R$, so its zero-th homotopy group is $\pi_0(R)^{\times}$ and its higher homotopy agrees with that of $R$.  Barring the issue of computing $\Pic^0(R)$, which in general is a difficult problem, this allows us to write down the homotopy groups of $B\Pic(R)$ in terms of those of $R$.  These groups are listed for a few common ring spectra in table~\ref{tablehtpy}.

\begin{table}[h]
\begin{centering}
\begin{tabular}{|c|cccccccc|}
\hline
Classifying Space & \multicolumn{8}{|c|}{Homotopy Groups}
\\
& $i=1$ & $2$ & $3$ & $4$ & $5$ & $6$ & $7$ & $8$ \\
\hline
$B\Pic(\SS)$ & $\Z$ & $\Z/2$ & $\Z/2$ & $\Z/2$ & $\Z/24$ & $0$ & $0$ & $\Z/2$ \\
\hline
$B\Pic(\HZ)$ & $\Z$ & $\Z/2$ & $0$ & $0$ & $0$ & $0$ & $0$ & $0$ \\
\hline
$B\Pic(K)$ & $\Z/2$ & $\Z/2$ & $0$ & $\Z$ & $0$ & $\Z$ & $0$ & $\Z$ \\
\hline
$B\Pic(MU)$ & $(\Z)$ & $\Z/2$ & $0$ & $\Z$ & $0$ & $\Z^2$ & $0$ & $\Z^3$ \\
\hline
\end{tabular}
\caption[The homotopy groups of the classifying spaces for haunts]{The homotopy groups $\pi_i(B\Pic(R))$ of the classifying spaces for $R$-haunts.} \label{tablehtpy}
\end{centering}
\end{table}

\nid The parenthetical group $\pi_1(B\Pic(MU))$ is conjectural.  We now describe three simple examples of specters and their potential homology invariants.  The first specter has no global homology type, the second has a uniquely determined global homology type, and the third has a global homology type that depends on the choice of trivialization.

\begin{example}
Consider the specter described in example~\ref{s1specter}: the base space is $S^1$, the haunt $L_n$ has monodromy $S^n$, and the specter $T$ is $S^n \times D^1 \ra D^1$ on one semicircle and the cone $C(S^n \sqcup S^0) \ra C(\pt \sqcup \pt) = D^1$ on the other semicircle.  The classifying map of $L_n$ represents $n \in \Z = \pi_1(B\Pic(\SS))$; in particular there is no global parametrized spectrum corresponding to $T$ and therefore no ordinary homotopy type.  Moreover, the map $B(\Pic(\SS)) \ra B(\Pic(\HZ))$ is an isomorphism on $\pi_1$, so the homology haunt $(L_n)_{\HZ}$ is still nontrivial; correspondingly, the homology specter $T_{\HZ}$ is not a parametrized $\HZ$-module and so $T$ does not have homology invariants.
\end{example}

\begin{example}
Next consider the specter from example~\ref{s3specter}: the base space is $S^3$, the haunt $L$ is the haunt determined by the equatorial transition function $- \sm_{S^2} (S^2 {\rtimes} S^2)$, and the specter $P$ is $S^0 \times D^3$ on one hemisphere and $C(S^2 {\rtimes} S^2)$ on the other.  On the one hand, the classifying map $S^3 \ra B\Pic(\SS)$ for $L$ is nontrivial, representing $1 \in \Z/2 = \pi_3(B\Pic(\SS))$; on the other hand, there are no nontrivial $\HZ$-haunts on $S^3$ and so $L_{\HZ}$ is trivializable.  Thus, even though the specter $P$ has no global homotopy type, it does have homology invariants.  Moreover, $\Pic(\HZ)$ has homotopy only in degrees 0 and 1, so the set of trivializations of $L_{\HZ}$, namely $[S^3,\Pic(\HZ)] \cong \Pic^0(\HZ) = \{\Sigma^n \HZ\}$, is as small as possible.  The homology invariants of $P$ are therefore uniquely determined up to degree shift, that is up to tensoring with $\pi_*(\Sigma^n \HZ)$.

Let us calculate the homology $H_*(P)$ of the specter $P$.  Roughly speaking, the spectrum $\HZ$ cannot see the difference between the transition functions $- \sm_{S^2} (S^2 {\rtimes} S^2)$ and $- \sm_{S^2} (S^2 \times S^2)$.  As a result, the parametrized $\HZ$-module $P_{\HZ}$ is equivalent to the one obtained by gluing together $\HZ \times D^3$ and $C(\Sigma^2 \HZ \times S^2)$ using the map $- \sm_{S^2} (S^2 \times S^2)$.  By desuspending the second hemisphere, this is in turn equivalent to the parametrized $\HZ$-module $(\HZ \times D^3) \cup_{S^2} (C(\HZ \times S^2))$.  This last $\HZ$-module is simply the reduced homology of $S^3$; thus $H_*(P)$ is $\Z$ in a single degree and zero in all other degrees.  We will formalize this sort of computation in a moment.
\end{example}

\begin{example}
Let $P$ be the parametrized spectrum $S^0 \times (S^2 \times S^1)$ over $S^2 \times S^1$.  This spectrum defines a specter for the trivial haunt $H$ and as such $P$ has homology invariants.  However, there are two distinct trivialization $\tau_0$ and $\tau_1$ of $H_{\HZ}$.  They yield respectively the parametrized $\HZ$-modules $\tau_0(P_{\HZ}) = \HZ \times (S^2 \times S^1)$ and $\tau_1(P_{\HZ})= \HZ {\rtimes} (S^2 \times S^1)$; this last module exhibits a mobius transformation of the fibre $\HZ$ along the $S^1$ factor of the base.  The homology groups of these two trivializations are $H_*^{\tau_0}(P) \cong \{\Z \: \Z \: \Z \: \Z\}$ and $H_*^{\tau_1}(P) \cong \{\Z/2 \: 0 \: \Z/2 \: 0\}$.  The difference between these two groups is characteristic of the ambiguity involved in $\HZ$-specters on non-simply connected base spaces.  This dependence on the trivialization indicates that $\HZ$-specter invariants are not naturally graded abelian groups but rather objects in a more subtle algebraic category.
\end{example}

Typically specters are constructed by specifying parametrized spectra over the open sets of a cover; these parametrized spectra do not agree on the intersections but instead are glued together by the transition functions of the haunt.  This explicit local presentation suggests a method for computing specter invariants in terms of the local homology invariants of the defining parametrized spectra:

\begin{prop}
Suppose $H$ is a haunt, on a connected space $X$, whose associated $R$-haunt $H_R$ admits a trivialization $\tau: H_R \xra{\sim} (X \times \Rmod)$.  Let $P$ be a specter for $H$ and let $P_R$ denote the associated $R$-specter.  Then there is a "Mayer-Vietoris" spectral sequence
\begin{equation} \nn
E^2_{pq} =  H_p(X; \underline{\pi_q(\tau(P_R))}) \Rightarrow R_{p+q}^{\tau}(P).
\end{equation}
Here $\underline{\pi_q(\tau(P_R))}$ denotes the cosheaf $U \mapsto \pi_q((\tau(P_R) |_U)/U)$.  

Suppose $\{U_i\}$ is a fixed contractible cover of $X$, and let $- \sm_{U_{ij}} \rho_{ij}$ be transition functions defining $H$.  Suppose the specter $P$ is presented by parametrized spectra $P_i$ on $U_i$ together with identifications $\gamma_{ij}: P_{i} |_{U_{ij}} \sm_{U_{ij}} \rho_{ij} \xra{\sim} P_j |_{U_{ij}}$.  Then the above spectral sequence has the form
\begin{equation} \nn
E^1_{pq} = \bigoplus_{i_1 < \ldots < i_p} R_q(P_{i_1}  |_{U_{i_1\ldots i_p}}) \Rightarrow R_{p+q}^{\tau}(P).
\end{equation}
If the trivialization is given by automorphisms $\tau_i: (U_i \times \Rmod) \xra{\sim} (U_i \times \Rmod)$ appropriately compatible with the $R$-haunt transition functions $(\rho_{ij})_R = \rho_{ij} \sm_{U_{ij}} R$, then the $d^1$ differential is given by the maps
\begin{align}
R_q(P_{i_1}  |_{U_{i_1\ldots i_p}}) &\xra{(\iota_k)_*} R_q(P_{i_1}  |_{U_{i_1\ldots \widehat{i_k} \ldots i_p}}) \nn\\
R_q(P_{i_1}  |_{U_{i_1\ldots i_p}}) &\xra{(\iota_1)_* \circ (\tau_{i_2})_*^{-1} \circ (\tau_{i_1})_*} R_q(P_{i_2}  |_{U_{i_2 \ldots i_p}}) \nn
\end{align}
Here $\iota_k: U_{i_1\ldots \widehat{i_k} \ldots i_p} \ra U_{i_1 \ldots i_p}$ and $\iota_1: U_{i_2 \ldots i_p} \ra U_{i_1 \ldots i_p}$ denote the inclusions.
\end{prop}

The first half of this proposition is just a statement, in ordinary parametrized homotopy theory, about the parametrized spectrum $\tau(P_R)$; it is not in itself particularly useful because one must expressly identify $\tau(P_R)$ in terms of the original specter $P$ in order to compute the cosheaf homology.  The second half is more explicit and addresses the situation that actually arises with twisted parametrized spectra.  In particular, the above $E^1$ term does not depend on the trivialization $\tau$ and can be immediately computed in any given case.

\section{Polarized Hilbert Manifolds and Semi-Infinite Spectra}  \label{sectionpolar} 

Thusfar our discussion has been purely topological: the category of spectra has a complicated space of automorphisms and it is natural to study bundles of categories of spectra and their associated sections, namely twisted parametrized spectra or "specters".  These bundles are, however, intimately connected to the geometry of infinite-dimensional manifolds, and homotopy-theoretic invariants of such manifolds often take the form of twisted parametrized spectra.  In the first part of this section, we describe the relevant geometry, namely real polarizations of real Hilbert bundles, and we show how a manifold equipped with this structure gives rise to a bundle of categories of spectra, that is to a haunt.  We also discuss a related structure, a "unitary" polarization and investigate the invariants of specters for haunts associated to unitary polarizations; the resulting description of these invariants provides an extensive generalization of the Cohen-Jones-Segal complex-oriented Floer invariants~\cite{cjs}.  In the second part of this section, we introduce a conjectural construction of the category of specters, for a given polarized bundle, in terms of parametrized semi-infinitely indexed spectra.  Specifically, instead of indexing spectra on finite-dimensional subspaces of a countably infinite-dimensional vector space, we introduce spectra indexed on the semi-infinite subspaces of a Hilbert space that are compatible with a fixed polarization.  A parametrized version of these semi-infinite spectra provides an explicit geometric viewpoint on the category of specters.

\subsection{The Homotopy Theory of Polarized Bundles} \label{htpypolar}

We describe four types of polarizations, namely real and complex on a real Hilbert space and real and complex on a complex Hilbert space, and discuss the homotopy types of the corresponding classifying spaces---it turns out that there is no distinction between the two notions of polarization on a complex Hilbert space.  We then describe polarized bundles and note that a real polarization on a real Hilbert bundle gives rise to a haunt.  We conclude by investigating the special class of "unitary" polarized bundles, namely real polarizations of a real bundle that lift to polarizations of a complex bundle.  In particular we show that under mild conditions, a specter for a unitary polarization admits $\HZ$-, $K$-, and $MU$-homology invariants.

\subsubsection{Polarizations of Hilbert Space} \label{polhilb}

A finite dimensional vector bundle on a space $X$ is classified by a map from $X$ to the classifying space $BO(n)$.  The topology of the classifying space is governed by the (non-trivial) topology of the orthogonal group $O(n)$ of automorphisms of $\R^n$.  By contrast, the orthogonal group $O(\HH)$ of Hilbert space is contractible~\cite{kuiper} and therefore the classifying space $BO(\HH)$ carries no topological information; indeed, all Hilbert bundles on a given space are isomorphic.  In particular, if $X$ is a Hilbert manifold, that is an infinite-dimensional manifold whose tangent bundle is a Hilbert bundle, then the tangent bundle of $X$ carries no information at all about the topology of $X$.  The situation is not as bad as it might seem, however, because many naturally occurring infinite-dimensional manifolds come equipped with a polarization.  This polarization is a reduction of the structure group of $X$ from the orthogonal group $O(\HH)$ to the so-called restricted orthogonal group $O_{\res}(\HH)$.  This latter group does have an interesting topology, and so we can recover information about such a polarized Hilbert manifold $X$ from its polarized tangent bundle.

There are various notions that go under the name "polarization" and we spend a moment describing and distinguishing them; references include~\cite{pressleysegal, segalwilson, cjs} but the reader is warned that the terminology and definitions in those papers disagree with one another and at points with our treatment.  We fix an infinite-dimensional separable real Hilbert space $\HH$.  Morally, a polarization of $\HH$ is an equivalence class of decompositions $V \oplus W$ of $\HH$ or $\HH_{\C}=\HH \otimes \C$ arising from an eigenvalue decomposition of an appropriate operator $J:\HH \ra \HH$.  The subspaces $V$ and $W$ are sums of collections of eigenspaces of $J$; that the polarization is an equivalence class of decompositions rather than a single decomposition reflects an ambiguity over whether to assign certain eigenspaces to $V$ or to $W$.  For example, if $J$ is a self-adjoint Fredholm operator, then the associated decompositions $V \oplus W$ are roughly those in which $V$ contains almost all the eigenspaces for negative eigenvalues of $J$ and $W$ contains almost all the eigenspaces for positive eigenvalues of $J$.  If on the other hand $J$ is a skew-adjoint Fredholm operator, then the decompositions are those in which $V$ contains almost all the eigenspaces for positive imaginary eigenvalues and $W$ contains almost all the eigenspaces for negative imaginary eigenvalues.  

In practice many polarizations arise from self- and skew-adjoint Fredholm operators as above, but we can simplify the definitions of polarizations if we restrict attention to self- and skew-adjoint orthogonal isomorphisms.  That is, suppose $J:\HH \ra \HH$ is a self-adjoint orthogonal isomorphism; in this case, $J^2=1$ and so $\HH$ is split into the $+1$ and $-1$ eigenspaces $V$ and $V^{\perp}$.  Of course, any  orthogonal decomposition arises as the eigenvalue decomposition of such an operator, and so decompositions $\HH = V \oplus V^{\perp}$ are in one-to-one correspondence with orthogonal operators $J$ with $J^2=1$.  An equivalence class of such decompositions defines a "real" polarization, as follows.

\begin{defn} \label{defsymppol}
A \emph{real polarization} on a real Hilbert space $\HH$ is a collection of orthogonal decompositions $\{\HH = V \oplus V^{\perp}\}$ satisfying the following conditions:
\begin{itemize}
\item both $V$ and $V^{\perp}$ are infinite dimensional,
\item for any two decompositions $V \oplus V^{\perp}$ and $W \oplus W^{\perp}$ in the collection, the projections $V \ra W$ and $V^{\perp} \ra W^{\perp}$ are Fredholm and the projections $V \ra W^{\perp}$ and $V^{\perp} \ra W$ are Hilbert-Schmidt,
\item any decomposition $W \oplus W^{\perp}$ satisfying the second property with respect to a decomposition $V \oplus V^{\perp}$ in the collection is in the collection.
\end{itemize}
\end{defn}

\nid From now on, whenever we mention a decomposition of a Hilbert space, we implicitly assume that both factors of the decomposition are infinite dimensional.  Corresponding to the above definition we have a restricted orthogonal group:

\begin{defn}
Let $V \oplus V^{\perp}$ be a fixed decomposition of the real Hilbert space $\HH$.  The \emph{real restricted orthogonal group} $O_{\res}^r(\HH)$ is the subgroup of the
orthogonal group $O(\HH)$ of operators $\phi$ such that $\phi(V) \oplus \phi(V^{\perp})$ is
in the same real polarization class as $V \oplus V^{\perp}$.
\end{defn}

\nid Note that this is \emph{not} the group that Pressley and Segal~\cite{pressleysegal} refer to as the restricted orthogonal group.  Indeed we will see later that it has a radically different homotopy type than their $O_{\res}(\HH)$.

\begin{note}
The space of real polarizations of a real Hilbert space is the quotient $O(\HH)/O_{\res}^r(\HH)$.
\end{note}

Now by contrast, suppose we had begun with a skew-adjoint orthogonal isomorphism $J:\HH \ra \HH$; in this case, $J^2 = -1$ and so $\HH_{\C}$ is decomposed into the $+i$ and $-i$ eigenspaces $W$ and $\overline{W}$.  Indeed, for any decomposition $\HH_{\C} = W \oplus \overline{W}$ of $\HH_{\C}$ into a subspace $W$ and its conjugate $\overline{W}$, the following three conditions are equivalent:

\begin{itemize}
\item the decomposition is orthogonal with respect to the Hermitian metric $\langle-,-\rangle$ on $\HH_{\C}$ extending the inner product on $\HH$,
\item the subspaces $W$ and $\overline{W}$ are isotropic with respect to the bilinear form $(a,b) = \langle a,\overline{b}\rangle$ on $\HH_{\C}$,
\item there is an orthogonal operator $J:\HH \ra \HH$ with $J^2=-1$ having $W$ and $\overline{W}$ as its $+i$ and $-i$ eigenspaces respectively.
\end{itemize}

\nid We refer to such decompositions as "orthogonal".  An equivalence class of these orthogonal decompositions defines a "complex" polarization:

\begin{defn}
A \emph{complex polarization} on a real Hilbert space $\HH$ is a collection of orthogonal decompositions $\{\HH_{\C} = W \oplus \overline{W}\}$ satisfying the same conditions as in definition~\ref{defsymppol}.  Given a fixed orthogonal decomposition $W \oplus \overline{W}$, the \emph{complex restricted orthogonal group} $O_{\res}^c(\HH)$ on the real Hilbert space $\HH$ is the subgroup of $O(\HH)$ of operators $\phi$ such that $\phi(W) \oplus \phi(\overline{W})$ is in the same complex polarization class as $W \oplus \overline{W}$.  
\end{defn}
\nid This complex restricted orthogonal group is what Pressley and Segal~\cite{pressleysegal} refer to as $O_{\res}(\HH)$.

\begin{note}
The space of complex polarizations of a real Hilbert space is the quotient $O(\HH)/O_{\res}^c(\HH)$.
\end{note}

There is yet another notion that goes under the name polarization.  Let $\HH$ now be a complex Hilbert space.  Suppose $J: \HH \ra \HH$ is a self-adjoint unitary operator; then $J^2=1$ and $\HH$ is decomposed into $+1$ and $-1$ eigenspaces.  Two such decompositions $V \oplus V^{\perp}$ and $W \oplus W^{\perp}$ are considered equivalent if, as in definition~\ref{defsymppol}, the projections $V \ra W$ and $V^{\perp} \ra W^{\perp}$ are Fredholm and the other two projections are Hilbert-Schmidt.  An equivalence class of these decompositions defines a polarization of the complex Hilbert space $\HH$; the group of unitary operators preserving such a polarization is called $U_{\res}(\HH)$ and the space of such polarizations is $U(\HH)/U_{\res}(\HH)$.  Similarly if $J$ is skew-adjoint unitary, then $J^2=-1$ and $\HH$ is decomposed into $+i$ and $-i$ eigenspaces.  Equivalence classes of these decompositions also give a notion of polarization, but because $\HH$ is complex there is a one-to-one correspondence between self- and skew-adjoint unitary operators and the two notions of polarization coincide.  This correspondence, which in a sense encodes the two-fold complex Bott periodicity, can obscure the distinction between the two notions in the real case.

We briefly discuss the homotopy types of these various spaces of polarizations.  As above, the real restricted orthogonal group $O_{\res}^r(\HH)$ of a real Hilbert space is the space of orthogonal operators $\phi: \HH \ra \HH$ such that the projections $\phi(V) \ra V$ and $\phi(V^{\perp}) \ra V^{\perp}$ are Fredholm and the projections $\phi(V) \ra V^{\perp}$ and $\phi(V^{\perp}) \ra V$ are Hilbert-Schmidt, for a fixed decomposition $V \oplus V^{\perp}$.  Suppose that we have chosen $\phi  |_V$ such that $\phi(V) \ra V$ is Fredholm and $\phi(V) \ra V^{\perp}$ is Hilbert-Schmidt.  Then the subspace $\phi(V^{\perp})$ is necessarily $\phi(V)^{\perp}$ and we can chose $\phi  |_{V^{\perp}}$ to be any orthogonal isomorphism $V^{\perp} \xra{\cong} \phi(V)^{\perp}$---the space of such choices is of course contractible.  Specifying $\phi  |_V$ amounts to choosing a Fredholm map $V \ra V$ and a Hilbert-Schmidt map $V \ra V^{\perp}$.  The space of Hilbert-Schmidt operators is contractible and so $O_{\res}^r(\HH)$ has the homotopy type of the space of Fredholm operators, namely $\Z \times BO$.  The associated space of polarizations $O(\HH)/O_{\res}^r(\HH)$ therefore has the homotopy type $B(\Z \times BO) \simeq U/O$.  By contrast, Pressley and Segal~\cite{pressleysegal} show that $O_{\res}^c(\HH)$ has the homotopy type $O/U \simeq \Omega O$ and so the space of polarizations $O(\HH)/O_{\res}^c(\HH)$ has the homotopy type $B(O/U) \simeq B(\Omega O) \simeq O$.  When $\HH$ is a complex Hilbert space, the restricted group $U_{\res}(\HH)$ has the homotopy type $\Z \times BU$ and the space of polarizations is $U(\HH)/U_{\res}(\HH) \simeq B(\Z \times BU) \simeq B(\Omega U) \simeq U$.  The homotopy types of these various spaces of polarizations are summarized in table~\ref{tablepolar}.

\begin{table}[!h]
\begin{centering}
\begin{tabular}{|c|c|c|}
\hline
Polarization Type & Structure Group & Classifying Space \\
\hline
real on $\HH_{\R}$ & $O_{\res}^r \simeq \Z \times BO$ & $B(\Z \times BO) \simeq U/O$ \\
\hline
complex on $\HH_{\R}$ & $O_{\res}^c \simeq \Omega O$ & $B(\Omega O) \simeq O$ \\
\hline
real on $\HH_{\C}$ & $U_{\res} \simeq \Z \times BU$ & $B(\Z \times BU) \simeq U$ \\
\hline
complex on $\HH_{\C}$ & $U_{\res} \simeq \Omega U$ & $B(\Omega U) \simeq U$ \\
\hline
\end{tabular}
\caption{The homotopy types of the classifying spaces for polarizations.} \label{tablepolar}
\end{centering}
\end{table}

We will be primarily concerned with real polarizations of real Hilbert spaces and unless otherwise indicated, "polarization" will refer to this notion.


\subsubsection{Real Polarizations and Haunts}  \label{symppol}

A priori a Hilbert bundle $E$ on a space $X$ is classified by a map $X \ra BO(\HH)$---of course this map contains no topological information.  To give a polarization of this bundle is to specify, continuously in $X$, a polarization on each fibre of $E$:

\begin{defn}
A polarization of the Hilbert bundle $E$ on the space $X$ is a reduction of the structure group of $E$ from $O(\HH)$ to $O_{\res}^r(\HH)$.  In other words it is a lift of the classifying map $X \ra BO(\HH)$ to a map $X \ra BO_{\res}^r(\HH)$.
\end{defn}

\nid As there is no harm in doing so, we usually think of a polarization on the Hilbert bundle $E$ simply as a map $X \ra BO_{\res}^r(\HH) \simeq B(\Z \times BO)$.  A polarization of a Hilbert manifold is simply a polarization of its (trivial) tangent bundle.  The fundamental link between the geometry of polarizations and twisted parametrized stable homotopy theory is the association
\begin{equation} \nn
\left\{ \textrm{\parbox{3.35cm}{\centering Polarized Hilbert \\ bundles on $X$}} \right\} \rightsquigarrow 
\left\{ \textrm{\parbox{1.4cm}{\centering Haunts \\ on $X$}} \right\}
\end{equation}

\nid This association is determined by composing the classifying map $X \ra B(\Z \times BO)$ of the polarized bundle with the deloop of the J-homomorphism $B(\Z \times BO) \xra{BJ} B(\Z \times BG)=B(\Z \times B\GL_1(\SS))$.

The basic philosophy behind this correspondence is that geometric structures on infinite dimensional manifolds are intimately connected with polarizations and that homotopy-theoretic information about these structures can be encoded in twisted parametrized spectra for the haunt associated to the polarization.  The specific nature of this connection will be the subject of future work with Mike Hopkins~\cite{douglashopkins} and Ciprian Manolescu~\cite{douglasmanolescu}.  Here we record a few illustrative examples of polarized manifolds and their associated haunts.

There are two widely utilized sources of polarized manifolds: the first is loop spaces of symplectic and almost complex manifolds, and the second is moduli spaces of connections in gauge theory---see for example~\cite{cjs}.  We discuss the first source of examples.
Given a symplectic manifold $M$, a choice of metric determines an almost complex structure $M \xra{c} BU(n)$ on the tangent bundle of $M$.  The loop of this classifying map, or indeed of the classifying map for any almost complex manifold, can be used to determine a polarization on $LM$:
\begin{equation} \nn
LM \xra{Lc} LBU(n) \ra \Omega BU(n) \ra \Omega BU \simeq U \ra U/O \simeq B(\Z \times BO)
\end{equation}
This polarization and its associated haunt can be highly nontrivial.

\begin{example} \label{poleg1}
Let $S^6$ have its usual almost complex structure.  The haunt associated to the resulting polarization of the loop space $LS^6$ is nontrivial.  Indeed, the classifying map
\begin{equation} \nn
LS^6 \ra U \ra U/O \ra B(\Z \times BG)
\end{equation}
for this haunt restricts on $S^5 \hookrightarrow LS^6$ to a generator of $\pi_5(B(\Z \times BG)) = \pi_3(\SS) = \Z/24$.  To see this, note that because $\pi_4(S^3) = \Z/2$, the Hurewicz map $\pi_6(BSU(3)) \ra H_6(BSU(3))$ is multiplication by 2.  Because the Euler characteristic of $S^6$ is 2, this implies that the almost complex structure $S^6 \ra BSU(3)$ is a generator of $\pi_6(BSU(3)) \cong \pi_6(BSU)$.  The loop $LS^6 \ra LBU$ of this almost complex structure therefore induces an isomorphism on $\pi_5$ and the claim follows:
\begin{equation} \nn
\xymatrix{
\pi_5(LS^6) \ar[r]  \ar@{=}[d] & \pi_5(LBU) \ar[r] \ar@{=}[d]
& \pi_5(U) \ar[r] \ar@{=}[d] & \pi_5(U/O) \ar[r] \ar@{=}[d] &
\pi_5(B(\Z \times BG)) \ar@{=}[d] \\
\Z \ar[r]^{\cong} & \Z \ar[r]^{\cong}
& \Z \ar[r]^{\cong} & \Z \ar@{->>}[r] & \Z/24
}
\end{equation}
A specter for the resulting canonical haunt on $LS^6$ will have no global homotopy type.  We will see in a moment though that any such specter has $\HZ$-, $K$-, and $MU$-homology invariants.
\end{example}

\subsubsection{Unitary Polarizations and Specter Invariants} \label{unitpol}

Many examples of polarized manifolds have the property that the polarization map $X \ra U/O$ factors through a map $X \ra U$; (this is true for instance of example~\ref{poleg1} above).  As we saw earlier, the stable unitary group $U$ classifies polarizations on complex Hilbert bundles $E$.  The projection map $U \ra U/O$ corresponds to viewing a polarization of $E$ as a real polarization of the underlying real Hilbert bundle $E_{\R}$; (similarly, the inclusion map $U \ra O$ corresponds to viewing the polarization of $E$ as a complex polarization of the underlying real Hilbert bundle $E_{\R}$).  For lack of better terminology, we say that a real polarization $X \ra U/O$ of a real Hilbert bundle is \emph{unitary} if it lifts to a polarization $X \ra U$ of a complex Hilbert bundle.

Haunts associated to unitary polarizations are much better behaved than arbitrary haunts in a sense we now describe, and as a result their specters have a much simplified invariant theory.  Suppose $X \ra U \ra U/O$ is the classifying map for a unitary polarization.  The associated haunt $H$ has a corresponding $\HZ$-haunt $H_{\HZ}$ classified by the composite
\begin{equation} \nn
X \ra U \ra U/O \ra B(\Z \times BG) \ra B(\Z \times B\Z/2) \simeq B\Z \times B^2\Z/2.
\end{equation}
The homology group $H^2(U;\Z/2)$ is zero, so this composite factors through $B\Z$.  Let $q \in \Z$ denote the smallest nonzero integer in the image of $H_1(X;\Z) \ra H_1(B\Z;\Z)$, and suppose $T$ is a specter for the haunt $H$ on $X$.  The haunt $H_{\HZ}$ is nontrivial; thus $T_{\HZ}$ does not have the form of a parametrized $\HZ$-module and so has no associated global $\HZ$-module, therefore no corresponding chain complex and no homology groups, per se.  
There is nevertheless a \emph{$\Z/q$-graded} chain complex associated to $T_{\HZ}$ and therefore $T$ has $\Z/q$-graded homology groups for invariants.  This fact nicely explains the idea (by now prevalent in the literature) that semi-infinite and Floer homology theories are naturally graded not by the integers but by finite cyclic groups.

If the unitary polarization $X \ra U \ra U/O$ is trivial on $H^1$, we have a more complete description of the corresponding haunt and its associated invariants:

\begin{prop}
Let $X$ be a space having the homotopy type of a finite CW complex and let $E$ be a Hilbert bundle on $X$ equipped with a unitary polarization $X \ra U \ra U/O$.  Suppose the induced map $H^1(U;\Z) \ra H^1(X;\Z)$ is zero.  Then any specter for the haunt $H$ associated to the polarized bundle $E$ admits global $\HZ$-, $K$-, and $MU$-homology invariants.
\end{prop}
\begin{proof}
Because of the $H^1$ condition, the classifying map for the polarization factors through $SU$, and because $X$ is homotopy finite, this map in turn factors through some $SU(n)$.  The haunt $H$ is therefore classified by a map $X \ra SU(n) \ra B(\Z \times BG)$.  Let $R$ denote one of the spectra $\HZ$, $K$, or $MU$.  The $R$-haunt $H_R$ is classified by the composition $X \ra SU(n) \ra B(\Z \times BG) = B\Pic(\SS) \ra B\Pic(R)$.  It is not of course the case that the map $B\Pic(\SS) \ra B\Pic(R)$ is null, but the composition $SU(n) \ra B\Pic(\SS) \ra B\Pic(R)$ will be null for the spectra $R$ in question.  This we can see by considering the map from the Atiyah-Hirzebruch spectral sequence for $B\Pic(\SS)^*(SU(n))$ to the Atiyah-Hirzebruch spectral sequence for $B\Pic(R)^*(SU(n))$.  If for example $R=K$ and $n=3$ the map of $E_2$ terms is as follows:

\begin{table}[!h]
\begin{centering}
{\footnotesize
\begin{tabular}{c|ccccccccc cccc |ccccccccc}
\cline{2-10} \cline{15-23}
0 & 0& 0 & 0 &0 & 0 &0 & 0 &0 & 0 & &&&& 0 & 0 &0 & 0 &0 & 0 & 0 & 0 & 0 
\\
1 & $\Z$ & 0 & 0 & $\Z$ & 0 & $\Z$ & 0 & 0 & $\Z$ & &&&& $\Z/2$ & 0 & 0 & $\Z/2$ & 0 & $\Z/2$ & 0 & 0 & $\Z/2$ 
\\
2 & $\Z/2$ & 0 & 0 & $\Z/2$ & 0 & $\Z/2$ & 0 & 0 & $\Z/2$ &&&&& $\Z/2$ & 0 & 0 & $\Z/2$ & 0 & $\Z/2$ & 0 & 0 & $\Z/2$ 
\\
3 & $\Z/2$ & 0 & 0 & \framebox{$\Z/2$} & 0 & $\Z/2$ & 0 & 0 & $\Z/2$ &&&&& 0 & 0 &0 & 0 &0 & 0 & 0 & 0 & 0 
\\
4 & $\Z/2$ & 0 & 0 & $\Z/2$ & 0 & $\Z/2$ & 0 & 0 & $\Z/2$ & \multicolumn{3}{c}{$\longrightarrow$}
&& $\Z$ & 0 & 0 & $\Z$ & 0 & $\Z$ & 0 & 0 & $\Z$
\\
5 & $\Z/24$ & 0 & 0 & $\Z/24$ & 0 & \framebox{$\Z/24$} & 0 & 0 & $\Z/24$ &&&&& 0 & 0 &0 & 0 &0 & 0 & 0 & 0 & 0
\\
6 & 0 & 0 &0 & 0 &0 & 0 & 0 & 0 & 0 &&&&& $\Z$ & 0 & 0 & $\Z$ & 0 & $\Z$ & 0 & 0 & $\Z$
\\
7 & 0& 0 & 0 &0 & 0 &0 & 0 &0 & 0 & &&&& 0 & 0 &0 & 0 &0 & 0 & 0 & 0 & 0 
\\
8 & $\Z/2$ & 0 & 0 & $\Z/2$ & 0 & $\Z/2$ & 0 & 0 & \framebox{$\Z/2$} &&&&& $\Z$ & 0 & 0 & $\Z$ & 0 & $\Z$ & 0 & 0 & \framebox{$\Z$}
\end{tabular}
}
\end{centering}
\end{table}

\nid Total degree zero terms are boxed.  The map is necessarily zero at $E_2$ in total degree zero and therefore $B\Pic(\SS)^0(SU(3)) \ra B\Pic(K)^0(SU(3))$ is zero.  The cases of $\HZ$ and $MU$ and of other $n$ are similar.  We therefore conclude that the $R$-haunt $X \ra B\Pic(R)$ is trivializable and so any corresponding specter has $R$-homology invariants.
\end{proof}

The proposition says that for a large class of polarized manifolds, the associated semi-infinite homotopy types, namely the specters, will have homology, $K$-theory, and complex bordism invariants.  This provides an explanation of and a substantial generalization of the remark in Cohen-Jones-Segal~\cite{cjs} that trivially polarized manifolds should have semi-infinite homology, $K$-theory, and complex bordism invariants.

\subsection{Semi-Infinitely Indexed Spectra} \label{siis}

In section~\ref{symppol} we saw that the classifying space for polarizations maps to the classifying space for haunts; any polarized bundle therefore has an associated haunt and a corresponding category of specters.  In this section we will sketch, using a notion of parametrized semi-infinitely indexed spectra, a conjectural realization of this category of specters in terms of the geometry of the polarized bundle.

We begin by defining (non-parametrized) semi-infinitely indexed spectra.  Classically a prespectrum $E$ is presented by specifying a space $E(\R^n)$ for each integer $n$ together with appropriate structure maps.  This notion naturally evolved (in work of May and company~\cite{mayquinnray, lewismaystein}) into that of coordinate free prespectra; a coordinate free prespectrum is presented by giving a space $E(V)$ for every finite dimensional subspace $V$ of a fixed countably-infinite dimensional inner product space $\R^{\infty}$, together with appropriate structure maps $\Sigma^{W-V} E(V) \ra E(W)$ for each inclusion $V \subset W$.  A prespectrum is a spectrum if the adjoint structure maps $E(V) \ra \Omega^{W-V} E(W)$ are homeomorphisms.  Such a spectrum has the property that if $V$ and $W$ have the same dimension, then $E(V)$ and $E(W)$ are homeomorphic, and moreover $E(V)$ varies continuously with $V$.  This can be seen by noting that for any compact family $\{V_t\}$ of finite-dimensional subspaces of $\R^{\infty}$, there is a finite-dimensional subspace $W$ such that $W$ contains all the subspaces of the family; the spaces $E(V_t)$ are therefore determined as $\Omega^{W-V_t} E(W)$, which is evidently a continuous family.  Later on, Elmendorf~\cite{elmendorf} included as part of the definition of a spectrum the requirement that the spaces $E(V)$ vary continuously in $V$.

The fundamental idea behind semi-infinitely indexed spectra is that the natural subspaces of a polarized Hilbert space $\HH$ are not the finite dimensional subspaces but the "negative energy" or "semi-infinite" subspaces, that is the subspaces $V$ such that $V \oplus V^{\perp}$ is a decomposition of $\HH$ in the given polarization class.  A semi-infinitely indexed spectrum is then roughly an assignment of a space $E(V)$ to each such semi-infinite $V$, together with appropriate structure maps.  It is not the case that given a compact family $\{V_t\}$ of semi-infinite subspaces, there exists a semi-infinite subspace $W$ containing all the $V_t$; indeed, there exist decompositions $V \oplus V^{\perp}$ and $W \oplus W^{\perp}$ in the same polarization class such that $V$ and $W$ span the whole Hilbert space $\HH$.  As such, it is important that we impose a continuity condition on the spaces $E(V)$---we do so roughly along the lines of~\cite{elmendorf} and we thank Mike Hopkins for bringing that reference to our attention.

\begin{defn} \label{X-spectrum}
Let $X$ be a space together with a distinguished subset $X^{(2)}$ of $X^2$ and a finite dimensional vector bundle $\gamma$ on $X^{(2)}$.  Define $X^{(3)}$ to be $\{(a,b,c) \in X^3  | (a,b),(a,c),(b,c) \in X^{(2)}\}$.  Let $p_{12}$, $p_{13}$, and $p_{23}$ denote the projections $X^{(3)} \ra X^{(2)}$ to the indicated factors, and similarly denote by $p_1$, $p_2$, and $p_3$ the projections $X^{(3)} \ra X$.  Suppose there is an identification $p_{13}^* \gamma = p_{12}^* \gamma \oplus p_{23}^* \gamma$.  Then an \emph{$X$-prespectrum} is a bundle $T$ on $X$ together with a map $\sigma: S^{\gamma} \sm_{X^{(2)}} p_1^* T \ra p_2^* T$ such that the diagram
\begin{equation} \nn
\xymatrix{
S^{p_{23}^* \gamma} \sm_{X^{(3)}} S^{p_{12}^* \gamma} \sm_{X^{(3)}} p_1^* T \ar@{=}[r] \ar[d]
& S^{p_{13}^* \gamma} \sm_{X^{(3)}} p_1^* T  \ar[d] \\
S^{p_{23}^* \gamma} \sm_{X^{(3)}} p_2^* T \ar[r] & p_3^* T
} 
\end{equation}
commutes.  This data forms an \emph{$X$-spectrum} if the adjoint of $\sigma$ is a homeomorphism.
\end{defn}

Suppose $\HH$ is equipped with a fixed polarization.  Let $Gr_{\res}(\HH)$ denote the grassmannian of decompositions $V \oplus V^{\perp}$ in the polarization class.  The space $Gr_{\res}^{(2)}(\HH)$ is the set of pairs of decompositions $(V \oplus V^{\perp}, W \oplus W^{\perp})$ such that $V \subset W$, and the vector bundle $\gamma$ is the orthogonal complement $V^{\perp W}$.

\begin{defn}
A \emph{semi-infinitely indexed (pre)spectrum}, or "semi-infinite (pre)spectrum" for short, is a $Gr_{\res}(\HH)$-(pre)spectrum.
\end{defn}

We bother with the abstract definition~\ref{X-spectrum} because it facilitates comparisons between $X$-spectra as $X$ varies.  In particular, fix a decomposition $\HH = \HH^- \oplus \HH^+$ in the polarization class and let $Gr(\HH^+)$ denote the grassmannian of finite dimensional subspaces of $\HH^+$.  We will refer to $Gr(\HH^+)$-spectra as Hilbert spectra.  Furthermore, denote by $\R^{\infty}$ a countably-infinite-dimensional dense subspace of $\HH^+$ and let $Gr(\R^{\infty})$ be the grassmannian of finite dimensional subspaces of $\R^{\infty}$.  We now have natural restriction maps
\begin{align}
Gr_{\res}(\HH){\textrm{-spectra}} & \ra Gr(\HH^+){\textrm{-spectra}} \nn \\
(V \subset \HH \rightsquigarrow E(V)) &\mapsto (W \subset \HH^+ \rightsquigarrow E(\HH^- \oplus W)) \nn \\
Gr(\HH^+){\textrm{-spectra}} & \ra Gr(\R^{\infty}){\textrm{-spectra}} \nn \\
(W \subset \HH^+ \rightsquigarrow F(W)) &\mapsto (U \subset \R^{\infty} \rightsquigarrow F(U)) \nn
\end{align}
A detailed treatment of the theory of semi-infinitely indexed spectra and their relationship to ordinary spectra will appear elsewhere---in particular one must construct a semi-infinite spectrification functor and semi-infinite sphere spectra (leading to a semi-infinite notion of stable weak equivalence) and one must build left adjoints to the two restriction maps above.  For now we leave as a conjecture the following:
\begin{conj} \label{sispectraconj}
There is a notion of weak equivalence of semi-infinite spectra and a notion of weak equivalence of Hilbert spectra such that the above restriction map from $Gr_{\res}(\HH)$-spectra to $Gr(\HH^+)$-spectra induces an equivalence of (simplicial) homotopy categories.  Similarly the restriction map from $Gr(\HH^+)$-spectra to $Gr(\R^{\infty})$-spectra induces an equivalence of homotopy categories.
\end{conj}

\twiddles

Parametrized semi-infinite spectra are no more difficult to define than semi-infinite spectra.  The idea of parametrized universes (of the countably-infinite variety) for spectra first appeared in Elmendorf~\cite{elmendorf}.  Though we were not aware of this reference during our development of parametrized semi-infinite spectra, it is a very clean presentation of the classical case and we follow it in spirit.  Elmendorf's real insight was not so much the use of parametrized universes, per se, as the realization that one could build a category of spectra on all (parametrized) universes at once and that this larger category was substantially better than the category of spectra indexed on a single universe.  Unfortunately, basic facts about countably-infinite universes that made this possible fail to be true of semi-infinite universes, and so we necessarily shy aware from this aspect of Elmendorf's treatment.

Let $E$ be a polarized Hilbert bundle on a space $X$.  The restricted grassmannian $Gr_{\res}(E)$ of this bundle is the space of semi-infinite subspaces $V$ of fibres $E_p$, $p \in  X$, of $E$; that is, the subspaces $V$ are such that $V \oplus V^{\perp}$ is a decomposition in the polarization class of $E_p$.  The space of pairs $Gr_{\res}^{(2)}(E)$ and the finite-dimensional bundle $\gamma$ are defined as before.  We immediately have

\begin{defn}
A \emph{parametrized semi-infinite spectrum} for the polarized bundle $E$ is a $Gr_{\res}(E)$-spectrum.
\end{defn}

Granting conjecture~\ref{sispectraconj}, the claim is quite simply that the category of parametrized semi-infinite spectra for the polarized Hilbert bundle $E$ on the space $X$ is homotopically equivalent to the category of specters (twisted parametrized spectra) for the haunt on $X$ associated to the polarized bundle.  Thereby, semi-infinite spectra provide a geometric realization of the homotopy-theoretic correspondence, via the classifying map \mbox{$B(\Z \times BO) \ra B(\Z \times BG)$}, between polarizations and haunts.

\subsection*{Acknowledgments}

As noted in the introduction, I am thoroughly indebted to Mike Hopkins, Bill Dwyer, and Jacob Lurie.  I would also like to thank Johann Sigurdsson for explaning various aspects of his work with Peter May on parametrized homotopy theory; Phil Hirschhorn and Mark Hovey for tolerating numerous questions about model categories; Charles Rezk for comments about homotopy limits of complete segal spaces; and Ciprian Manolescu, Andre Henriques, Stefan Schwede, Jeff Smith, and Stephan Stolz for lending an ear and thoughtfully responding to various stages of this project.

\bibliography{tpshtv3}
\bibliographystyle{plain}

\end{document}